\documentclass[11pt,a4paper,reqno]{amsart}
\usepackage{amssymb,amsmath,amsthm}
\usepackage[utf8]{inputenc}
\usepackage[T1]{fontenc}
\usepackage{enumerate}
\usepackage[all]{xy}
\usepackage{fullpage}
\usepackage{comment}
\usepackage{array}
\usepackage{longtable}
\usepackage{stmaryrd}
\usepackage{mathrsfs} 
\usepackage{xcolor}
\usepackage{mathtools}
\usepackage[foot]{amsaddr}

\def\deg{\mathrm{deg}}

\def\int{\mathrm{Int}}

\def\dim{\mathrm{dim}}
\def\ba{\mathbf{a}}

\usepackage{hyperref}
\hypersetup{
    colorlinks=true,
    linkcolor=blue,
    filecolor=red,  
    citecolor=green,
    urlcolor=cyan,
    pdftitle={On Covering Radii in function fields},
    pdfpagemode=FullScreen,
    }

\usepackage{cleveref}
\crefformat{section}{\S#2#1#3}
\crefformat{subsection}{\S#2#1#3}
\crefformat{subsubsection}{\S#2#1#3}
\usepackage{enumitem}
\usepackage{tikz}

\newtheorem{theorem}{Theorem}[section]
\newtheorem{lemma}[theorem]{Lemma}
\newtheorem{corollary}[theorem]{Corollary}
\newtheorem{proposition}[theorem]{Proposition}

\newtheorem{conjecture}[theorem]{Conjecture}

\theoremstyle{definition}
\newtheorem{definition}[theorem]{Definition}

\theoremstyle{remark}
\newtheorem{remark}[theorem]{Remark}

\makeatletter
\newcommand{\subalign}[1]{%
  \vcenter{%
    \Let@ \restore@math@cr \default@tag
    \baselineskip\fontdimen10 \scriptfont\tw@
    \advance\baselineskip\fontdimen12 \scriptfont\tw@
    \lineskip\thr@@\fontdimen8 \scriptfont\thr@@
    \lineskiplimit\lineskip
    \ialign{\hfil$\m@th\scriptstyle##$&$\m@th\scriptstyle{}##$\hfil\crcr
      #1\crcr
    }%
  }%
}
\makeatother

\numberwithin{equation}{section}

\title{On Covering Radii in function fields}
\author[1]{Noy Soffer Aranov}
\address{Department of Mathematics, Technion, Technion City, Haifa, Israel 3200003}
\email{noyso@campus.technion.ac.il}

\begin{document}

\maketitle
\begin{abstract}
    In this paper, we shall discuss topics in geometry of numbers in the function field setting, such as covering radii. We find a closed form for covering radii with respect to convex bodies, which will lead to a proof of the function field analogue of Woods' conjecture in this setting. Then, we will prove a function field analogue of Minkowski's conjecture about the multiplicative covering radius. To do this, we shall prove a function field analogue of Solan's result that every diagonal orbit intersects the set of well rounded lattices. This implies that the Gruber-Mordell spectrum in function field is trivial in every dimension.
    

\end{abstract}

\section{Introduction}
Let $d\geq 2$, let $G=\operatorname{SL}_d(\mathbb{R})$, let $\Gamma=\operatorname{SL}_d(\mathbb{Z})$, and let $X_d=G/\Gamma$. Then $X_d$ can be identified with the space of unimodular lattices via the identification
$$g\Gamma\mapsto g\mathbb{Z}^d$$
Then, $G$ acts naturally on $X_d$ by $g\cdot h\mathbb{Z}_d=gh\mathbb{Z}^d$. In this paper, we shall discuss several function field analogues of fundamental theorems in geometry of numbers. This builds upon and extends the results of \cite{Mah}. One such theorem is Minkowski's convex body theorem.
\begin{theorem}
    \label{MCBT}
    Let $\mathcal{C}$ be a convex set that is symmetric around the origin such that $\operatorname{Vol}(\mathcal{C})>2^d$. Then, for every $\Lambda\in X_d$, $\Lambda\cap \mathcal{C}\neq \{0\}$.
\end{theorem}
Given $\mathcal{C}$ convex and symmetric around the origin and a lattice $\Lambda\in X_d$, one can ask how well Theorem \ref{MCBT} can be improved for this particular $\mathcal{C}$. In other words, given a convex symmetric set $\mathcal{C}$, we can ask what is the largest $r>0$ such that $\Lambda\cap r\mathcal{C}=\{0\}$. This question has been well studied when $\mathcal{C}$ is a symmetric box. Explicitly, given a lattice $\Lambda\in X_d$ and a symmetric box of the form 
$$\mathcal{B}=[-a_1,a_1]\times ,\dots, \times [-a_d,a_d],$$
we say that $\mathcal{B}$ is admissible for $\Lambda$ if $\mathcal{B}\cap \Lambda=\{0\}$. We define the Mordell constant of $\Lambda$ as
$$\kappa(\Lambda)=\frac{1}{2^d}\sup_{\mathcal{B}\text{ admissible}}\operatorname{Vol}(\mathcal{B}).$$
Define the Gruber-Mordell spectrum as
$$\mathbf{MG}_d=\{\kappa(\Lambda):\Lambda\in X_d\}.$$
Theorem \ref{MCBT} implies that $\kappa(\Lambda)\leq 1$ for every $\Lambda\in X_d$. Moreover, it is easy to see that the box $B=(-1,1)^d$ is admissible for the lattice $\mathbb{Z}^d$. Hence, $\kappa(\mathbb{Z}^d)=1$. Furthermore, $\kappa$ is invariant under the group of diagonal matrices with determinant $1$, 
$$A=\Bigg\{\begin{pmatrix}
    a_1&&\\
    &\ddots&\\
    &&a_d
\end{pmatrix}:\prod_{i=1}^d\vert a_i\vert=1\Bigg\}.$$
Since $A$ acts ergodically on $X_d$, then, $\kappa$ is constant almost everywhere. In addition, by \cite[Corollary 1.3]{SW15}, $\kappa(\Lambda)=1$ for almost every $\Lambda\in X_d$.

Often times, it will be interesting to understand whether a ball contains at least $k$ linearly independent lattice points. A way to describe this is through successive minima. For $r>0$ and $x\in \mathbb{R}^d$, we denote 
$$B(x,r)=\{y\in \mathbb{R}^d:\Vert x-y\Vert_2<r\}.$$
For $i=1,\dots,d$ and a lattice $\Lambda=g\mathbb{Z}^d\subseteq \mathbb{R}^d$, we define the $i$-th successive minima of $\Lambda$ as
$$\lambda_i(\Lambda)=\inf\{r>0:\Lambda\cap B(0,r) \text{ contains }r \text{ linearly independent vectors}\}.$$
We denote $\det(\Lambda)=\vert \det(g)\vert$.  
\begin{theorem}[Minkowski's 2nd Theorem]
\label{Mink2ndR}
Let $\Lambda=g\mathbb{Z}^d$ be a lattice. Then,
$$\frac{2^d}{d!\operatorname{Vol}(B(0,1))}\vert \det(g)\vert\leq \prod_{i=1}^d\lambda_i(\Lambda)\leq \frac{2^d}{\operatorname{Vol}(B(0,1))}\vert \det(g)\vert.$$
\end{theorem}
Theorem \ref{Mink2ndR} can also be interpreted geometrically. The determinant of a lattice is bounded by above an below by two universal constants times the volume of the parallelopiped formed by a set of linearly independent vectors, whose lengths correspond to the successive minima of a lattice. 

Another example of a function which provides geometric information about a lattice is the covering radius. For $\Lambda\in X_d$ and a compact convex body $\mathcal{C}\subseteq \mathbb{R}^d$, we define the covering radius of $\Lambda$ with respect to $\mathcal{C}$, $\operatorname{CovRad}_{\mathcal{C}}(\Lambda)$ as the smallest $r\geq 0$ such that $\Lambda+r\mathcal{C}=\mathbb{R}^d$. Since $\mathcal{C}$ is compact, then the function $\operatorname{CovRad}_{\mathcal{C}}:X_d\rightarrow \mathbb{R}$ is a proper function. In particular, it is not bounded by above. In \cite{Rog}, Rogers proved that 
$$\inf_{x\in X_d}\operatorname{CovRad}_{\mathcal{C}}(x)\leq d^{\log_2\log d+O(1)}$$
In \cite{ORW}, Ordentlich, Regev and Weiss improved Rogers' bound on $\operatorname{CovRad}_{\mathcal{C}}$. 
\begin{theorem}{\cite[Theorem 1.1]{ORW}}
    \label{ORW}
    For any compact convex body $\mathcal{C}$, 
    $$\inf_{\Lambda\in X_d}\operatorname{CovRad}_{\mathcal{C}}(\Lambda)=O(d^2).$$
\end{theorem}
The functions $\Theta_{d,\mathcal{C}}=\inf_{\Lambda\in X_d}\operatorname{CovRad}_{\mathcal{C}}(\Lambda)$ and $\Theta_d=\inf_{\Lambda\in X_d}\operatorname{CovRad}_{B(0,1)}(\Lambda)$ have been studied extensively \cite{Rog,CS,GL,ORW}. In general, it is interesting to understand covering radii of lattices with respect to functions. Given a function $F:\mathbb{R}^d\rightarrow \mathbb{R}^+$ with $F(0)=0$ and a lattice $\Lambda\in X_d$, we define $\operatorname{CovRad}_F(\Lambda)$ to be the smallest $r>0$ such that for every $R>r$,
\begin{equation}
    \Lambda+\{\mathbf{v}\in \mathbb{R}^d:F(\mathbf{v})<R\}=\mathbb{R}^d.
\end{equation}
For example if $F(\mathbf{v})=\Vert \mathbf{v}\Vert_2$, then, $\operatorname{CovRad}_F=\operatorname{CovRad}_{B(0,1)}$. An interesting covering radius is the Minkowski function, $\mu=\operatorname{CovRad}_N$, where $N(\mathbf{v})=\prod_{i=1}^d\left|v_i\right|$. It is clear that $\mu$ is invariant under the group $A$, and therefore, due to ergodicity of the $A$-action \cite{Mau}, $\mu$ is constant almost everywhere. Furthermore due to upper semicontinuity, the generic value of $\mu$ is its lower bound, which is $0$, as Shapira showed in \cite{S09}.
\begin{conjecture}[Minkowski]
\label{MinkR}
For every $\Lambda\in X_d$, we have that
\begin{enumerate}
    \item \label{minkUppR} $\mu(\Lambda)\leq 2^{-d}=\mu(\mathbb{Z}^d)$ and
    \item \label{minkUnique}$\mu(\Lambda)=2^{-d}$ if and only if $\Lambda\in A\mathbb{Z}^d$.
\end{enumerate}
\end{conjecture}
Conjecture \ref{MinkR} has been proved for $d\leq 10$ \cite{Mink,Rem,Dys,Sku,W,HRS7,HRS8,KR9,KR10,Sol}. One method to prove Conjecture \ref{MinkR}(\ref{minkUppR}) is to use $A$ invariance of $\mu$ to prove that the conclusion of Conjecture \ref{MinkR}(\ref{minkUppR}) holds for a subset of $X_d$ which intersects every $A$ orbit.
\begin{theorem}
    Assume that there exists a set $Y\subseteq X_d$ such that
    \begin{enumerate}
        \item \label{MinkConjY} For every $\Lambda\in Y$, $\mu(\Lambda)\leq q^{-d}$.
        \item \label{AorbIntersectR} For every $\Lambda\in X_d$, $A\Lambda\cap Y\neq \emptyset$.
    \end{enumerate}
    Then, every $\Lambda\in X_d$ satisfies that $\mu(\Lambda)\leq 2^{-d}$.
\end{theorem}
This method was utilized in \cite{McM,LSW,Sol} when $Y$ is the set of well rounded lattices.
\begin{definition}
    For $\Lambda\in X_d$, we say that $\Lambda$ is well rounded if $\lambda_1(\Lambda)=\dots=\Lambda_d(\Lambda)$. We denote by $\mathcal{WR}_d$ the set of well rounded lattices in $X_d.$
\end{definition}
Solan \cite{Sol} proved that for every $d\geq 2$, the set of well rounded lattices satisfy (\ref{AorbIntersectR}). Thus, to conclude the proof of Conjecture \ref{MinkR})(\ref{minkUppR}), it suffices to prove that the well rounded lattices satisfy (\ref{MinkConjY}). One method to prove that a subset $Y\subseteq X_d$ satisfies condition (\ref{MinkConjY}) is to prove that $\operatorname{CovRad}_{\Vert \cdot\Vert_2}$ is bounded from above by $\frac{\sqrt{d}}{2}$.
\begin{lemma}
\label{normCov->Mink}
    Assume that $\Lambda\in X_d$ satisfies $\operatorname{CovRad}_{\Vert \cdot\Vert_2}(\Lambda)\leq \frac{\sqrt{d}}{2}$. Then, $\mu(\Lambda)\leq q^{-d}$.
\end{lemma}
Woods \cite{W} conjectured that the set of well rounded lattices satisfy the conditions of Lemma \ref{normCov->Mink}. Regev, Shapira, and Weiss \cite{RSW} proved that Woods' conjecture is false for $d\geq 30$. In \cite{SW16}, Shapira and Weiss conjectured that the stable lattices, which are lattices whose subgroups all have covolume at least $1$, satisfy the conditions of Lemma \ref{normCov->Mink}. If this is indeed true, this would conclude the proof of Conjecture \ref{MinkR}(\ref{minkUppR}), due to \cite{Sol}.
 
In this paper, we shall discuss geometry of numbers in function field. In particular, we shall provide a closed form of the function field analogue of $\operatorname{CovRad}_{\mathcal{C}}$ and prove a function field analogue of Theorem \ref{ORW}. Furthermore, we shall also obtain a complete description of the function field analogue of the function $\mu$ in dimension $2$ and the function $\kappa$ in every dimension, as well as prove a function field analogue of Conjecture \ref{MinkR}(\ref{minkUppR}). Several of these results will stem from a decomposition theorem, which improves the function field analogue of Theorem \ref{Mink2ndR}.
\subsection{Definitions}
\label{subsec:Def}
Let $d\geq 2$, let $p$ be a prime, $q$ be a power of $p$, and let $\mathcal{R}=\mathbb{F}_q[x]$ be the ring of polynomials over $\mathbb{F}_q$. Let $\mathcal{K}=\mathbb{F}_q(x)$ be the field of rational functions over $\mathbb{F}_q$. We define an absolute value on $\mathcal{R}$ by $\vert f\vert=q^{\deg(f)}$ and extend it to an absolute value on $\mathcal{K}$ by 
$$\left|\frac{f}{g}\right|=\begin{cases}
    q^{\deg(f)-\deg(g)}&f\neq 0\\
    0&f=0
\end{cases}.$$ 
Then the topological completion of $\mathcal{K}$ with respect to the metric $d(f,g)=\vert f-g\vert$ is the field of Laurent series $\tilde{\mathcal{K}}$ defined by
$$\tilde{\mathcal{K}}\:= \mathbb{F}_q\left(\left(x^{-1}\right)\right)=\Bigg\{\sum_{n=-N}^{\infty} a_n x^{-n}:a_n\in \mathbb{F}_q,N\in \mathbb{Z}\Bigg\}.$$
Given $d\geq 2$, we define the norm on $\tilde{\mathcal{K}}^d$ by $\Vert \mathbf{v}\Vert=\max_{i=1,\dots,d}\vert v_i\vert.$ 
\begin{lemma}{Ultrametric Inequality}
\label{lem:UM}
    \begin{enumerate}
            \item For all $\alpha,\beta\in \mathbb{F}_q((x^{-1}))$, we have $\vert \alpha+\beta\vert\leq \max\{\vert \alpha\vert,\vert \beta\vert\}$.
            \item For all $\mathbf{u},\mathbf{v}\in \mathbb{F}_q((x^{-1}))^n$, we have $\Vert\mathbf{u}+\mathbf{v}\Vert\leq \max\{\Vert \mathbf{u}\Vert,\Vert\mathbf{v}\Vert\}$. 
        \end{enumerate}        
\end{lemma}
Let $\mathcal{O}$ be the maximal compact order of $\tilde{\mathcal{K}}$, that is
$$\mathcal{O}=\mathbb{F}_q\left[\left[x^{-1}\right]\right]=\{f\in \tilde{\mathcal{K}}:\vert f\vert\leq 1\}.$$
Denote the group of units in $\mathcal{O}$ by $\mathbf{U}$, that is 
$$\mathbf{U}\:= \{f\in \tilde{\mathcal{K}}:\vert f\vert=1\}=\Bigg\{\sum_{n=0}^{\infty}a_nx^{-n}:a_n \in \mathbb{F}_q, a_0\in \mathbb{F}_q^*\Bigg\}=\mathcal{O}^*.$$
We can view $\tilde{\mathcal{K}}^*$ as the direct product $\tilde{\mathcal{K}}^*\cong \mathbb{Z}\times \mathbf{U}$ in the following way:
$$f\mapsto \left(\log_q\vert f\vert,\frac{f}{x^{\log_q\vert f\vert}}\right).$$
\begin{definition}
\label{def:rhopi}
Define the functions $\rho(f)=\log_q\vert f\vert$ and $\pi(f)=\frac{f}{x^{\log_q\vert f\vert}}$. 
By abuse of notation, we write $\rho(\mathbf{v})=(\rho(v_1),\dots, \rho(v_d))$ and similarly $\pi(\mathbf{v})=(\pi(v_1),\dots, \pi(v_d))$ for vectors $\mathbf{v}\in \tilde{\mathcal{K}}^d$. Similarly, for $g\in \operatorname{GL}_d(\tilde{\mathcal{K}})$ we define $(\rho(g))_{ij}=\rho(g_{ij})$ and $(\pi(g))_{ij}=\pi(g_{ij})$.
\end{definition}

Let $G=\operatorname{SL}_d(\tilde{\mathcal{K}})$ be the group of invertible $d\times d$ matrices over $\tilde{\mathcal{K}}$ with $\vert \det(g)\vert=1$. Let $\Gamma=\operatorname{SL}_d(\mathcal{R})<G$ be the group of invertible $d\times d$ matrices with entries in $\mathcal{R}$ with $\vert \det(g)\vert=1$. Let $\mathcal{L}_d=G/\Gamma$. Then, $\mathcal{L}_d$ is identified with the space of unimodular lattices in $\tilde{\mathcal{K}}^d$ under the identification
$$g\Gamma\mapsto g\mathcal{R}^d.$$
Let $A$ be the group of diagonal matrices in $G$ with determinant $1$. 
\begin{definition}
\label{length}
Given a lattice $\Lambda=g\Gamma\subseteq \tilde{\mathcal{K}}^d$, we define the length of the shortest non-zero vector in $\Lambda$ as
\begin{equation}
    \ell(\Lambda)=\min\big\{\Vert \mathbf{v}\Vert:\mathbf{v}\in g\Gamma\setminus\{0\}\big\},
\end{equation}
where $\Vert (v_1,\dots, v_d)^t\Vert=\max_i\vert v_i\vert$. 
\end{definition}
In $\mathcal{L}_d$, Mahler's compactness criterion gives a necessary and sufficient condition for compactness (see \cite{Cas59} for the real case, \cite[Theorem 3.3]{Gho} for the function field case, and \cite[Theorem 1.1]{KST} for general $S$-adic fields). 
\begin{theorem}[Mahler's Compactness Criterion]
\label{Mahler}
A closed set of lattices $Y\subseteq \mathcal{L}_d$ is compact if and only if there exists $\varepsilon>0$ such that $\inf_{\Lambda\in Y}\ell(\Lambda)>\varepsilon$. 
\end{theorem}
We shall now state the function field analogue of Theorem \ref{MCBT}. To do so, we define the Haar measure on $\tilde{\mathcal{K}}^d$ as the unique left invariant measure $m$ on $\tilde{\mathcal{K}}^d$, such that $m(\mathcal{O}^d)=1$. 
Mahler \cite{Mah} proved that every convex set in $\tilde{\mathcal{K}}^d$ is a linear image of $\mathcal{O}^d$.
\begin{theorem}[\cite{Mah}]
\label{ConvexSet}
    A set $\mathcal{C}\subseteq \tilde{\mathcal{K}}^d$ is an open compact $\mathcal{O}$-module, if and only if there exists $h\in \operatorname{GL}_d(\tilde{\mathcal{K}})$ such that $\mathcal{C}=h\mathcal{O}^d$.
    In this case we say that $\mathcal{C}$ is a convex body. We denote $m(\mathcal{C})=\operatorname{Vol}(\mathcal{C})=\vert \det(h)\vert$.
\end{theorem}
Notice that due to the definition of a convex body, any convex body must contain $0$. We can now state the function field analogue of Theorem \ref{MCBT}, which is a corollary of Serre Duality (see for example \cite{Cla}).
\begin{theorem}
\label{ClaMCBT}
    Let $\mathcal{C}\subseteq \tilde{\mathcal{K}}^d$ be a convex set and let $\Lambda\in \mathcal{L}_d$ be a lattice. If $m(\mathcal{C})\geq q^{-(d-1)}$, then $\mathcal{C}\cap \Lambda\neq \{0\}$. 
\end{theorem}
In function field, we have the following version of Theorem \ref{Mink2ndR}. 
\begin{theorem}{\cite{Mah}}
\label{Mink2nd}
For any lattice $\Lambda=g\mathcal{R}^d\subseteq \tilde{\mathcal{K}}^d$ and $i=1,\dots d$, we define the $i$-th successive minima of $\Lambda$ as
$$\lambda_i(\Lambda)=\min\{r>0: \text{ there exist }i \text{ linearly independent vectors in } \Lambda \text{ of norm }\leq r\}.$$
Then,
\begin{equation}
    \prod_{i=1}^d \lambda_i(\Lambda)=\vert \det(g)\vert.
\end{equation}
\end{theorem}
Theorem \ref{Mink2nd} can be generalized for general convex bodies (for the proof, see for example \cite[Corollary 5.3]{BK}). 
\begin{theorem}
    \label{Mink2ndConv}
    For a lattice $\Lambda=g\mathcal{R}^d\subseteq \tilde{\mathcal{K}}^d$, a convex body $\mathcal{C}$, and $i=1\dots d$, define
    $$\lambda_{i,\mathcal{C}}(\Lambda)=\min\{r>0:\text{ there exist }i \text{ linearly independent vectors }\Lambda \text{ in } r\mathcal{C}\}.$$
    Then, 
    \begin{equation}
        \prod_{i=1}^d\lambda_{i,\mathcal{C}}(\Lambda)=\frac{\vert \det(g)\vert}{\operatorname{Vol}(\mathcal{C})}.
    \end{equation}
\end{theorem}
Another well studied notion in parametric geometry of numbers is orthogonality \cite{KST,RW,AB24}.
\begin{definition}
    We say that $\mathbf{v}_1,\dots,\mathbf{v}_m\in \tilde{\mathcal{K}}^d$ are orthogonal if 
    \begin{equation}
        \Vert \mathbf{v}_1\wedge \dots\wedge\mathbf{v}_m\Vert=\prod_{i=1}^d\Vert \mathbf{v}_i\Vert.
    \end{equation}
    Equivalently, for any $\alpha_1,\dots,\alpha_m\in \tilde{\mathcal{K}}$, we have 
    $$\Vert\alpha_1\mathbf{v}_1+\dots+\alpha_m\mathbf{v}_m\Vert=\max_{i=1,\dots,m}\Vert \alpha_i\mathbf{v}_i\Vert,$$
    that is the equality case of Lemma \ref{lem:UM} holds. 
\end{definition}
By Theorem \ref{Mink2nd}, if $\mathbf{v}_1,\dots,\mathbf{v}_d$ satisfy $\Vert \mathbf{v}_i\Vert=\lambda_i(\Lambda)$, then, $\mathbf{v}_1,\dots,\mathbf{v}_d$ are orthogonal.  
\subsection{Main Results}
\subsubsection{Bases of Lattices}
\label{subsec:bases}
\begin{definition}
    We say that a set of linearly independent vectors $\mathbf{v}^{(1)},\dots,\mathbf{v}^{(d)}\in \Lambda$ is a set of successive minima vectors for $\Lambda$ if $\big\Vert \mathbf{v}^{(i)}\big\Vert=\lambda_i(\Lambda)$ for each $i=1,\dots,d$.  
\end{definition}
\begin{lemma}
\label{shortSpan}
Let $\Lambda\subseteq \tilde{\mathcal{K}}^d$ be a lattice and let $\mathbf{v}^{(1)},\dots,\mathbf{v}^{(d)}\in \Lambda$ be a set of successive minima vectors for $\Lambda$. Then, $\Lambda$ is spanned over $\mathcal{R}$ by $\mathbf{v}^{(1)},\dots,\mathbf{v}^{(d)}$.
\end{lemma}
\begin{remark}
    Due to Lemma \ref{shortSpan}, from now on, we shall call a set of successive minima vectors for $\Lambda$ a basis of successive minima for $\Lambda$.
\end{remark}
\begin{remark}
Lemma \ref{shortSpan} does not hold over $\mathbb{R}$. For example, take the lattice
$$L=\{\mathbf{v}\in \mathbb{R}^d:\forall 1\leq i,j\leq d, v_i\equiv v_j\mod \mathbb{Z}\}=2\mathbb{Z}^d+\mathbb{Z}\begin{pmatrix}
1\\
\vdots\\
1
\end{pmatrix}\in X_d.$$
If $d\geq 5$, then, $2\mathbf{e}_1,\dots, 2\mathbf{e}_d$ are the shortest linearly independent vectors of $L$. On the other hand, the vector $\begin{pmatrix}
1\\
\vdots\\
1
\end{pmatrix}\in L$ is not spanned over $\mathbb{Z}$ by $2\mathbf{e}_1,\dots, 2\mathbf{e}_d$.  
\end{remark}
From Theorem \ref{Mink2nd} and Lemma \ref{shortSpan}, we deduce that any set of successive minima for a lattice are an orthogonal $\mathcal{R}$ basis for the lattice. Moreover, all orthogonal $\mathcal{R}$ bases for a lattice are a basis of successive minima. 
\begin{theorem}
\label{thm:UniqOrtBase}
    Let $\Lambda\subseteq \tilde{\mathcal{K}}^d$ be a lattice. Then, every orthogonal $\mathcal{R}$ basis of $\Lambda$ is a basis of successive minima of $\Lambda$.
\end{theorem}
\subsubsection{Covering Radii with Respect to Convex Bodies}
\label{subsubsec:CovRad}
In this paper, we shall discuss several covering radii.
\begin{definition}
Let $d\geq 2$. 
\begin{enumerate}
    \item Given a convex body $\mathcal{C}\subseteq \tilde{\mathcal{K}}^d$ and a lattice $\Lambda\subseteq \tilde{\mathcal{K}}^d$, we define $\operatorname{CovRad}_{\mathcal{C}}(\Lambda)$ to be the smallest $r>0$ such that
\begin{equation}
\label{eqn:ConvCovRadDef}
    \Lambda+r\mathcal{C}=\tilde{\mathcal{K}}^d.
\end{equation}
We define the norm $\Vert \cdot \Vert_{\mathcal{C}}$ on $\tilde{\mathcal{K}}^d$ by $\Vert \mathbf{v}\Vert_{\mathcal{C}}=\inf\{r\geq 0:\mathbf{v}\in r\mathcal{C}\}$.
    \item In general, given a function $F:\tilde{\mathcal{K}}^d\rightarrow \mathbb{R}_+$ with $F(0)=0$, we can define $\operatorname{CovRad}_F(\Lambda)$ to be the smallest $r>0$ such that for every $R>r$, 
    \begin{equation}
    \label{eqn:CCovRad}
        \Lambda+\{\mathbf{v}\in \tilde{\mathcal{K}}^d:F(\mathbf{v})<R\}=\tilde{\mathcal{K}}^d.
    \end{equation}
\end{enumerate}
\end{definition}
\begin{remark}
    We can reinterpret $\operatorname{CovRad}_{\mathcal{C}}$ and $\lambda_{i,\mathcal{C}}$ using the terminology $\Vert \cdot \Vert_{\mathcal{C}}$.
    \begin{enumerate}
        \item We have that $\operatorname{CovRad}_{\mathcal{C}}=\operatorname{CovRad}_{\Vert \cdot \Vert_{\mathcal{C}}}$
        \item For $i=1\dots d$, $\lambda_{i,\mathcal{C}}(\Lambda)$ is the $i$-th successive minima of $\Lambda$ when measured with respect to the norm $\Vert \cdot \Vert_{\mathcal{C}}$. 
    \end{enumerate}
\end{remark}
We shall use the following stronger version of Theorem \ref{Mink2nd} to provide a closed form of $\operatorname{CovRad}_{\mathcal{C}}$ for a convex body $\mathcal{C}$. This can be viewed as a direct consequence of Grothendieck's classification of vector bundles over $\mathbb{P}^1$. We shall provide a self contained proof of this result.  
\begin{theorem}
\label{ughForm}
Let $g\in \operatorname{GL}_d(\tilde{\mathcal{K}})$ and let $\Lambda=g\mathcal{R}^d$. Then, there exist $u\in \operatorname{SL}_d(\mathcal{O})$ and $h\in \operatorname{SL}_d(\mathcal{R})$ such that
\begin{equation}
    ugh=\operatorname{diag}\Big\{x^{\log_q\lambda_1(\Lambda)},\dots ,x^{\log_q\lambda_d(\Lambda)}\Big\}.
\end{equation}
\end{theorem}
From Theorem \ref{ughForm}, we shall conclude the following claim regarding the covering radius with respect to $\mathcal{O}^d$.
\begin{theorem}
\label{NormCovRad} $\operatorname{CovRad}_{\mathcal{O}^d}(\Lambda)=q^{-1}\lambda_d(\Lambda)$.
\end{theorem}
\begin{corollary}
\label{ConvCovRad}
    Let $\mathcal{C}=h\mathcal{O}^d$ be a convex body and let $\Lambda=g\mathcal{R}^d$ be a lattice. Then,
    $$\operatorname{CovRad}_{\mathcal{C}}(\Lambda)=\frac{1}{q}\lambda_d\left(h^{-1}\Lambda\right)=\frac{1}{q}\lambda_d\left(h^{-1}g\mathcal{R}^d
    \right).$$
\end{corollary}
\begin{proof}
    Notice that (\ref{eqn:ConvCovRadDef}) holds if and only if 
    \begin{equation}
        h^{-1}g\mathcal{R}^d+r\mathcal{O}^d=\tilde{\mathcal{K}}^d.
    \end{equation}
    Hence, by Theorem \ref{NormCovRad}, $\operatorname{CovRad}_{\mathcal{C}}(\Lambda)=\operatorname{CovRad}_{\mathcal{O}^d}\left(h^{-1}g\mathcal{R}^d\right)=\frac{1}{q}\lambda_d\left(h^{-1}g\mathcal{R}^d\right)$.
\end{proof}
Furthermore, Corollary \ref{ConvCovRad} along with Theorem \ref{Mink2nd} enable us to prove a function field analogue of the main theorem in \cite{Rog,ORW}.
\begin{theorem}
\label{infConvCovRad}
For any convex body $\mathcal{C}$ with $\operatorname{Vol}(\mathcal{C})=1$, $\inf_{\Lambda\in \mathcal{L}_d}\operatorname{CovRad}_{\mathcal{C}}(\Lambda)=\frac{1}{q}$. 
\end{theorem}
\subsubsection{The Minkowski Function}
One application of Theorem \ref{NormCovRad} is proving the function field analogue of Conjecture \ref{MinkR}(\ref{minkUppR}). We first define the Minkowski function in function fields. Define $N:\tilde{\mathcal{K}}^d\rightarrow \mathbb{R}_+$ by $N(\mathbf{v})=\prod_{i=1}^d\left|v_i\right|$. For $\Lambda\in \mathcal{L}_d$, we define the Minkowski function by $\mu(\Lambda)=\operatorname{CovRad}_N(\Lambda)$. Explicitly, we can define $\mu$ in the following manner
$$\mu(\Lambda)=\sup_{\mathbf{v}\in \tilde{\mathcal{K}}^d}\inf_{\mathbf{u}\in \Lambda}N(\mathbf{v}-\mathbf{u}).$$
For $d\geq 2$, we define the Minkowski Spectrum by
$$\mathcal{S}_d=\{\mu(\Lambda):\Lambda\in \mathcal{L}_d\}.$$
It is easy to see that $\mu$ is invariant under the group of diagonal matrices with determinant of absolute value $1$, which we denote by $A$. Moreover, $\mu(\mathcal{R}^d)=q^{-d}$. Furthermore, we can bound $\mu$ by above in every dimension, which is a function field analogue of Conjecture \ref{MinkR}(\ref{minkUppR}). 
\begin{theorem}
    \label{MinkConj}
    For every $\Lambda\in \mathcal{L}_d$, $\mu(\Lambda)\leq q^{-d}$.
\end{theorem}
We note that due to \cite[Theorem 1.9]{SA}, for every $d\geq 2$, there are infinitely many  lattices in distinct $A$-orbits $A\Lambda\neq A\mathcal{R}^d$ such that $\mu(\Lambda)=q^{-d}$. Hence, the function field analogue of Conjecture \ref{MinkR}(\ref{minkUnique}) does not hold. Furthermore, a corollary of the main theorem from \cite{Agg} is that in dimension $2$, $\mu$ is a trivial function.  
\begin{theorem}
    \label{Mink2}
    $\mathcal{S}_2=\{q^{-2}\}$.
\end{theorem}
Similarly to the real case, a good strategy towards proving Theorem \ref{MinkConj} is to use the fact that $\mu$ is $A$-invariant in order to bound $\mu$ on a nice set of lattices instead of on all of $\mathcal{L}_d$. 
\begin{theorem}
\label{proofStrat}
    Assume that $Y\subseteq \mathcal{L}_d$ satisfies the following properties:
    \begin{enumerate}
        \item \label{Md}For every $\Lambda\in Y$, $\mu(\Lambda)\leq q^{-d}$.
        \item \label{Od}For every $\Lambda\in \mathcal{L}_d$, $A\Lambda\cap Y\neq \emptyset$.
    \end{enumerate}
    Then, Theorem \ref{MinkConj} holds in dimension $d$.
\end{theorem}
In this paper, we shall prove that when $Y$ is the set of well rounded lattices, the conditions of Theorem \ref{proofStrat} hold. This will lead to a proof of Theorem \ref{MinkConj}. 
\begin{definition}
    We say that $\Lambda\in \mathcal{L}_d$ is well rounded if $\Lambda$ contains $d$ linearly independent vectors of norm $\ell(\Lambda)$. We denote the set of well rounded lattices by $\mathcal{WR}_d$. In particular, due to Theorem \ref{Mink2nd}, $\lambda_j(\Lambda)=1$ for every $j=1\dots d$. 
\end{definition}
Theorem \ref{NormCovRad} implies the following function field analogue of Woods' conjecture. This is in stark contrast to the real case, in which Woods' conjecture does not hold for dimensions $d\geq 30$ \cite{RSW}. This is a consequence of Theorem 
\ref{Mink2nd}. 
\begin{remark}
    In \cite{SW16}, Shapira and Weiss conjectured that Woods' conjecture holds when replacing well rounded lattices with stable lattices. In our setting, well rounded and stable lattices coincide. 
\end{remark}
\begin{corollary}
    \label{Woods}
    For every $\Lambda\in \mathcal{WR}_d$, we have $\operatorname{CovRad}_{\mathcal{O}^d}(\Lambda)=q^{-1}$.
\end{corollary}
From here we obtain that every $\Lambda\in \mathcal{WR}_d$ must satisfy the conclusion of Theorem \ref{MinkConj}.
\begin{theorem}
\label{mu(WR)}
    For every $\Lambda\in \mathcal{WR}_d$, we have $\mu(\Lambda)\leq q^{-d}$.
\end{theorem}
\begin{proof}
    We first notice that $r\mathcal{O}^d=\{\mathbf{v}\in \tilde{\mathcal{K}}^d:\Vert \mathbf{v}\Vert\leq r\}$. Hence, $\operatorname{CovRad}_{\mathcal{O}^d}=\operatorname{CovRad}_{\Vert \cdot\Vert}$. Therefore, by Corollary \ref{Woods}, for every $\Lambda\in \mathcal{WR}_d$,
    \begin{equation}
        \Lambda+\{\mathbf{v}\in \tilde{\mathcal{K}}^d:\Vert \mathbf{v}\Vert\leq q^{-1}\}=\tilde{\mathcal{K}}^d.
    \end{equation}
    If $\Vert \mathbf{v}\Vert\leq q^{-1}$, then $N(\mathbf{v})=\prod_{i=1}^d\left|v_i\right|\leq q^{-d}$. Therefore,
    \begin{equation}
        \Lambda+\{\mathbf{v}\in \tilde{\mathcal{K}}^d:N(\mathbf{v})\leq q^{-d}\}=\tilde{\mathcal{K}}^d.
    \end{equation}
    Hence, for every $\Lambda\in \mathcal{WR}_d$, $\mu(\Lambda)\leq q^{-d}$.
\end{proof}
To conclude the proof of Theorem \ref{MinkConj}, we prove that the set of well rounded lattices must intersect every $A$ orbit. This can be viewed as a function field analogue of the main result of \cite{Sol}.
\begin{theorem}
\label{AxcapWR}
For every $\Lambda\in \mathcal{L}_d$, $A\Lambda\cap \mathcal{WR}_d\neq \emptyset$.
\end{theorem}
\begin{proof}[Proof of Theorem \ref{MinkConj}]
Let $\Lambda\in \mathcal{L}_d$. Then, by Theorem \ref{AxcapWR}, there exists $\mathbf{a}\in A$ such that $\mathbf{a}\Lambda\in \mathcal{WR}_d$. Hence, by Theorem \ref{mu(WR)} and $A$ invariance of $\mu$, $\mu(\Lambda)=\mu(\mathbf{a}\Lambda)\leq q^{-d}$. 
\end{proof}
Due to Theorem 1.9 in \cite{SA}, Theorem \ref{MinkConj}, Theorem \ref{Mink2}, and computer calculations, we propose the following conjecture.
\begin{conjecture}
\label{minktriv}
    For every $d\geq 2$, $\mathcal{S}_d=\{q^{-d}\}$.
\end{conjecture}

\subsubsection{The Mordell Function}
In addition, we shall prove that the function field analogue of the Mordell function is trivial. Given a box of the form 
$$\mathcal{B}=B(0,r_1)\times ,\dots,\times B(0,r_d)\subseteq \tilde{\mathcal{K}}^d,$$
and a lattice $\Lambda\in \mathcal{L}_d$, we say that $\mathcal{B}$ is admissible for $\Lambda$ if $\mathcal{B}\cap \Lambda=\{0\}$. We define the Mordell function $\kappa:\mathcal{L}_d\rightarrow \mathbb{R}^+$ by
$$\kappa(\Lambda)=\sup_{\mathcal{B}\text{ admissible}}\operatorname{Vol}(\mathcal{B}).$$
We define the Gruber Mordell spectrum by
$$\mathbf{MG}_d=\{\kappa(\Lambda):\Lambda\in \mathcal{L}_d\}.$$
Then, by Theorem \ref{ClaMCBT}, $\kappa(\Lambda)\leq q^{-d}$. Like its real analogue $\kappa$ is $A$-invariant and lower semicontinuous. Therefore, it is constant almost everywhere and its generic value is its upper bound. We shall prove that in $\mathcal{L}_d$, $\kappa$ is a trivial function which a direct consequence of Theorem \ref{ClaMCBT} and Theorem \ref{AxcapWR}.
\begin{theorem}
    \label{GM}
    For every $d\geq 2$, $\mathbf{MG}_d=\big\{q^{-d}\big\}$.
\end{theorem}
We note that Theorem \ref{GM} can be viewed as an easy corollary of Lemma 6.2 in \cite{BK} and Theorem \ref{ClaMCBT} but we shall include an alternative proof, since this proof uses dynamical tools and an improvement to Dirichlet's Theorem.  

\subsection{Acknowledgements}
I would like to thank Moshe Teutsch for help with writing the Sage code used in \cref{subsec:compMu}, Erez Nesharim for useful discussions regarding Theorem \ref{Mink2}. I would also like to thank Omri Nisan Solan for discussions leading to the proof of Theorem \ref{AxcapWR}. I would also like to thank Angelot Behajaina for reading a preliminary version of this paper and providing helpful comments. Most of all, I would like to thank Uri Shapira for introducing me to these questions and for countless discussions over the years about geometry of numbers, which enabled this paper to be possible. This work is supported by the ERC grant "Dynamics on Homogeneous Spaces" (no. 754475). I would also like to thank the anonymous referee for his comments, which helped make this paper clearer.
\section{Bases and Covering Radii with Respect to Convex Bodies}
In this section we shall prove the results in \cref{subsubsec:CovRad}. We shall first prove Lemma \ref{shortSpan} and deduce Theorem \ref{NormCovRad} and Theorem \ref{infConvCovRad} from Theorem \ref{ughForm}.
\begin{proof}[Proof of Lemma \ref{shortSpan}]
Assume towards a contradiction that there exists some $\mathbf{u}\in \Lambda$ which is not spanned over $\mathcal{R}$ by $\mathbf{v}^{(1)},\dots, \mathbf{v}^{(d)}$. Since $\mathbf{v}^{(1)},\dots, \mathbf{v}^{(d)}$ span $\tilde{\mathcal{K}}^d$ over $\tilde{\mathcal{K}}$, then there exist $n_i\in \mathcal{R}$ and $\vert \alpha_i\vert<1$, which are not all $0$, such that
\begin{equation}
    \mathbf{u}=\sum_{i=1}^dn_i\mathbf{v}^{(i)}+\sum_{i=1}^d\alpha_i\mathbf{v}^{(i)}.
\end{equation}
Notice that $\sum_{i=1}^d n_i\mathbf{v}^{(i)}\in \Lambda$ and thus 
$$0\neq \sum_{i=1}^d \alpha_i\mathbf{v}^{(i)}=\mathbf{u}-\sum_{i=1}^d n_i\mathbf{v}^{(i)}\in \Lambda.$$ 
Moreover, since $\mathbf{u}\notin \operatorname{span}_{\mathcal{R}}\big\{\mathbf{v}^{(1)},\dots, \mathbf{v}^{(d)}\big\}$, then $\sum_{i=1}^d \alpha_i\mathbf{v}^{(i)}\notin \operatorname{span}_{\mathcal{R}}\big\{\mathbf{v}^{(1)},\dots, \mathbf{v}^{(d)}\big\}$. On the one hand, by Lemma \ref{lem:UM}, 
\begin{equation}
    \Bigg\Vert \sum_{i=1}^d\alpha_i\mathbf{v}^{(i)}\Bigg\Vert\leq \max_i\vert \alpha_i\vert \Vert \mathbf{v}^{(i)}\Vert\leq \frac{1}{q}\lambda_d(\Lambda).
\end{equation}
On the other hand, since $\bigg\Vert \sum_{i=1}^d\alpha_i\mathbf{v}^{(i)}\bigg\Vert<\lambda_d(\Lambda)$ and $\sum_{i=1}^d \alpha_i\mathbf{v}^{(i)}\in \Lambda$, then, 
$$\sum_{i=1}^d\alpha_i\mathbf{v}^{(i)}\in \operatorname{span}_{\tilde{\mathcal{K}}}\big\{\mathbf{v}^{(1)},\dots,\mathbf{v}^{(d-1)}\big\}.$$ 
Hence $\alpha_d=0$. On the other hand, 
\begin{equation}
    \Bigg\Vert \sum_{i=1}^{d-1}\alpha_i\mathbf{v}^{(i)}\Bigg\Vert\leq \max_{i=1,\dots,d-1}\vert \alpha_i\vert\cdot \big\Vert \mathbf{v}^{(i)}\big\Vert\leq \frac{1}{q}\lambda_{d-1}(\Lambda)<\lambda_{d-1}(\Lambda).    
\end{equation}
Hence, $\sum_{i=1}^{d-1}\alpha_i\mathbf{v}^{(i)}\in \operatorname{span}_{\tilde{\mathcal{K}}}\big\{\mathbf{v}^{(1)},\dots,\mathbf{v}^{(d-2)}\big\}$ and thus $\alpha_{d-1}=0$. We can continue in such a manner to obtain that $\alpha_i=0$ for every $i=1,\dots,d$, which is a contradiction to the assumption that $(\alpha_1,\dots,\alpha_d)\neq (0,\dots, 0)$.
\end{proof}
\begin{proof}[Proof of Theorem \ref{thm:UniqOrtBase}]
    Let $\{\mathbf{v}_1,\dots,\mathbf{v}_d\}$ be an orthogonal basis for $\Lambda$ satisfying 
    $$\Vert \mathbf{v}_1\Vert\leq \dots\leq \Vert \mathbf{v}_d\Vert.$$ 
    Assume towards a contradiction that $\{\mathbf{v}_1,\dots,\mathbf{v}_d\}$ is not a basis of successive minima for $\Lambda$. Then, there exists a minimal $i_0$ such that $\Vert \mathbf{v}_{i_0}\Vert>\lambda_{i_0}(\Lambda)$. Hence, in order for 
$$\prod_{i=1}^d\Vert \mathbf{v}_i\Vert=\Vert \mathbf{v}_1\wedge \dots\wedge \mathbf{v}_d\Vert=\det(\Lambda),$$ there must exist some $i_1>i_0$ such that $\Vert \mathbf{v}_{i_1}\Vert<\lambda_{i_1}(\Lambda)$. On the other hand, by the definition of the successive minima, $\operatorname{span}(B(0,\Vert \mathbf{v}_{i_1}\Vert)\cap \Lambda)\subseteq \operatorname{span}(B(0,\lambda_{i_1}(\Lambda))$ has dimension at most $i_1-1$. Thus, $\dim\operatorname{span}\{\mathbf{v}_1,\dots,\mathbf{v}_{i_1}\}\leq i_1-1$, which contradicts the assumption that $\mathbf{v}_1,\dots,\mathbf{v}_{i_1}$ are orthogonal. Therefore, for every $i=1,\dots,d$, we have $\Vert \mathbf{v}_i\Vert=\lambda_i(\Lambda)$.
\end{proof}
\subsection{Proofs Assuming Theorem \ref{ughForm}}
We can obtain Theorem \ref{NormCovRad} as a corollary of Theorem \ref{ughForm}.
\begin{proof}[Proof of Theorem \ref{NormCovRad}]
Since $\operatorname{CovRad}$ is invariant under $\operatorname{SL}_d(\mathcal{O})$, then, by Theorem \ref{ughForm}, it suffices to prove Theorem \ref{NormCovRad} for $\Lambda=\operatorname{diag}\Big\{x^{\log_q\lambda_1(\Lambda)},\dots, x^{\log_q\lambda_d(\Lambda)}\Big\}\mathcal{R}^d$. 

Thus for any $\mathbf{v}\in \tilde{\mathcal{K}}^d$ there exist $n_i\in \mathcal{R}$ and $\vert \alpha_i\vert\leq q^{-1}$ such that $$\mathbf{v}=\sum_{i=1}^d(n_i+\alpha_i)x^{\log_q\lambda_i(\Lambda)}\mathbf{e}^{(i)}.$$ 
Let $\mathbf{w}=\sum_{i=1}^dn_i\mathbf{e}^{(i)}\in \Lambda$. Then, by the Lemma \ref{lem:UM}
\begin{equation}
\label{eqn:vDis}
    \Vert \mathbf{v}-\mathbf{w}\Vert=\bigg\Vert \sum_{i=1}^dx^{\log_q\lambda_i(\Lambda)}\alpha_i\mathbf{e}^{(i)}\bigg\Vert\leq \max_{i=1,\dots, d}\lambda_i(\Lambda)\vert \alpha_i\vert\cdot \big\Vert \mathbf{e}^{(i)}\big\Vert\leq \frac{1}{q}\lambda_d(\Lambda).
\end{equation}
Hence, (\ref{eqn:vDis}) implies that $\operatorname{CovRad}(\Lambda)\leq \frac{1}{q}\lambda_d(\Lambda)$. Let $\mathbf{v}=x^{\log_q\lambda_d(\Lambda)-1}\mathbf{e}^{(d)}\in \Lambda$. We shall show that 
\begin{equation}
    \min_{\mathbf{w}\in \Lambda}\Vert \mathbf{v}-\mathbf{w}\Vert=\Vert \mathbf{v}\Vert=\frac{1}{q}\lambda_d(\Lambda).
\end{equation}
If $\mathbf{w}\in \Lambda$, then, there exists $a_1,\dots,a_d\in \mathcal{R}$ such that $\mathbf{w}=\sum_{i=1}^dx^{\log_q\lambda_i(\Lambda)}a_i\mathbf{e}^{(i)}$. Hence
\begin{equation}
\label{eqn:v-wNorm}
    \Vert \mathbf{v}-\mathbf{w}\Vert=\max\Bigg\{\max_{i=1,\dots,d-1}\vert a_i\vert\lambda_i(\Lambda),\left|a_d-\frac{1}{x}\right|\lambda_d(\Lambda)\Bigg\}.
\end{equation}
Notice that (\ref{eqn:v-wNorm}) is minimal when $a_1=\dots=a_d=0$, and then, $\Vert \mathbf{v}-\mathbf{w}\Vert=\frac{1}{q}\lambda_d(\Lambda)$. Hence, 
\begin{equation}
\label{eqn:CovRad>=}
    \operatorname{CovRad}_{\mathcal{O}^d}(\Lambda)\geq\Vert \mathbf{v}-\mathbf{w}\Vert\geq \Vert \mathbf{v}\Vert=\frac{1}{q}\lambda_d(\Lambda),
\end{equation}
so that $\operatorname{CovRad}_{\mathcal{O}^d}(\Lambda)=\frac{1}{q}\lambda_d(\Lambda)$.
\end{proof}
We shall now obtain Corollary \ref{Woods} as a consequence of Theorem \ref{Mink2nd} and Theorem \ref{NormCovRad}.
\begin{proof}[Proof of Corollary \ref{Woods}]
Let $\Lambda\in \mathcal{WR}_d$. Then, by Theorem \ref{Mink2nd}, $\lambda_d(\Lambda)=1$ and thus, Theorem \ref{NormCovRad} implies that $\operatorname{CovRad}_{\mathcal{O}^d}(\Lambda)=\frac{1}{q}$.  
\end{proof}
Now Corollary \ref{ConvCovRad} and Theorem \ref{Mink2nd} imply Theorem \ref{infConvCovRad}.
\begin{proof}[Proof of Theorem \ref{infConvCovRad}]
    Let $\mathcal{C}=h\mathcal{O}^d$ be with $\operatorname{Vol}(\mathcal{C})=1$. Then, $\vert \det(h)\vert=1$. Thus, Corollary \ref{ConvCovRad} implies that for any $g\in \operatorname{SL}_d(\tilde{\mathcal{K}})$, $\operatorname{CovRad}_{\mathcal{C}}(g\mathcal{R}^d)=\frac{1}{q}\lambda_d(h^{-1}g\mathcal{R}^d)$. Since $h^{-1}g\mathcal{R}^d$ is unimodular, then, by Theorem \ref{Mink2nd}, $\lambda_d(h^{-1}g\mathcal{R}^d)\geq 1$. Furthermore, by taking $g=h$, we obtain that $\lambda_d(h^{-1}g\mathcal{R}^d)=1$, so that $\min_{\Lambda\in \mathcal{L}_d}\operatorname{CovRad}_{\mathcal{C}}(\Lambda)=\frac{1}{q}$. 
\end{proof}
\subsection{Proof of Theorem \ref{ughForm}}
Theorem \ref{ughForm} follows from Grothendieck's classification of vector bundles. We shall provide a self contained proof of Theorem \ref{ughForm} which does not use algebraic geometry.
\begin{proof}[Proof of Theorem \ref{ughForm}]
Let $\Lambda=g\mathcal{R}^d\subseteq \tilde{\mathcal{K}}^d$. By Lemma \ref{shortSpan}, there exists a basis of successive minima for $\Lambda$, which we denote by $\mathbf{v}^{(1)},\dots,\mathbf{v}^{(d)}$. By Theorem \ref{Mink2nd}, $\det(\Lambda)=\prod_{i=1}^d\big\Vert \mathbf{v}^{(i)}\big\Vert$. Let $g'$ be the matrix whose columns are $\mathbf{v}^{(1)},\dots,\mathbf{v}^{(d)}$. Then, $g'\mathcal{R}^d=g\mathcal{R}^d$, and therefore, $h=g^{-1}g'\in \operatorname{SL}_d(\mathcal{R})$. Let $u'=g'\operatorname{diag}\Big\{x^{-\log_q\lambda_1(\Lambda)},\dots, x^{-\log_q\lambda_d(\Lambda)}\Big\}$ be the matrix with columns $x^{-\log_q\lambda_1(\Lambda)}\mathbf{v}^{(1)},\dots, x^{-\log_q\lambda_d(\Lambda)}\mathbf{v}^{(d)}$. Then, $\vert \det(u')\vert=1$ and $u'\in M_{d\times d}(\mathcal{O})$. Hence, $u'\in \operatorname{SL}_d(\mathcal{O})$. Let $u=(u')^{-1}\in \operatorname{SL}_d(\mathcal{O})$. Thus,
\begin{equation}
    ugh=\operatorname{diag}\Big\{x^{\log_q(\lambda_1(\Lambda)},\dots, x^{\log_q\lambda_d(\Lambda)}\Big\}h^{-1}g^{-1}gh=\operatorname{diag}\Big\{x^{\log_q\lambda_1(\Lambda)},\dots,x^{\log_q\lambda_d(\Lambda)}\Big\}.
\end{equation}
As a consequence, 
$$u\Lambda=\operatorname{diag}\Bigg\{x^{\log_q\lambda_1(\Lambda)},\dots ,x^{\log_q\lambda_d(\Lambda)}\Bigg\}\mathcal{R}^d.$$
\end{proof}
\section{Well Rounded Lattices and the Minkowski Function}
\subsection{Proof of Theorem \ref{AxcapWR}}
\label{subsec:muUpp}
To prove Theorem \ref{AxcapWR}, we shall first show that the well rounded lattices must have a very particular form.
\begin{lemma}
\label{WRForm}
$\mathcal{WR}_d=\operatorname{SL}_d(\mathcal{O})\mathcal{R}^d$.
\end{lemma}
\begin{proof}
Let $\Lambda=g\mathcal{R}^d\in \mathcal{WR}$. By Theorem \ref{Mink2nd}, $\lambda_1(\Lambda)=1$. By Lemma \ref{shortSpan}, there exists a basis for $\Lambda$ over $\mathcal{R}$ of the form $\big\{\mathbf{v}^{(1)},\dots,\mathbf{v}^{(d)}\big\}$ such that $\big\Vert \mathbf{v}^{(i)}\big\Vert=1$ for every $i=1,\dots, d$. Let $g\in \operatorname{SL}_d(\tilde{\mathcal{K}})$ be defined by $\mathbf{v}^{(i)}=g\mathbf{e}^{(i)}$. Then, $\Lambda=g\mathcal{R}^d$. Since $\vert (g\mathbf{e}^{(i)})_j\vert\leq 1$ for every $i,j$, then $g$ must have entries in $\mathcal{O}$ and thus, $g\in \operatorname{SL}_d(\mathcal{O})$ so that $\Lambda\in \operatorname{SL}_d(\mathcal{O})\mathcal{R}^d$. 

On the other hand, if $g\in \operatorname{SL}_d(\mathcal{O})$, then, $g$ is an isometry. Hence, for every $\mathbf{v}\in \mathcal{R}^d$, we have
$\Vert g\mathbf{v}\Vert=\Vert \mathbf{v}\Vert\geq 1$. On the other hand, for every $i=1,\dots, d$, $\Vert g\mathbf{e}^{(i)}\Vert=1$ and therefore, $g\mathbf{e}^{(1)},\dots,g\mathbf{e}^{(d)}$ form a set of successive minima vectors for $g\mathcal{R}^d$. Hence, by Lemma \ref{shortSpan}, $g\mathcal{R}^d$ is spanned over $\mathcal{R}$ by $g\mathbf{e}^{(1)},\dots,g\mathbf{e}^{(d)}$, so that $g\mathcal{R}^d$ is well rounded.
\end{proof}
In order to prove Theorem \ref{AxcapWR}, we shall use the methods of \cite{Sol}. We shall first introduce some terminology which will enable us to reinterpret Theorem \ref{AxcapWR} in the language of flags.
\subsubsection{Wedge Products and Flags}
We shall need some terminology pertaining to wedge products. Let $\mathbf{e}^{(1)},\dots,\mathbf{e}^{(d)}$ be the standard basis of $\tilde{\mathcal{K}}^d$. 
\begin{definition}
For $1\leq k\leq d$, define $\bigwedge^k\tilde{\mathcal{K}}^d$ be the $k$-th exterior product of $\tilde{\mathcal{K}}^d$. We define the standard basis for $\bigwedge^k\tilde{\mathcal{K}}^d$ by $\mathbf{e}_J=\mathbf{e}^{(j_1)}\wedge \dots\mathbf{e}^{(j_k)}$ where $J=\{j_1<\dots <j_k\}\subseteq \{1,\dots,d\}$. Given $\mathbf{v}\in \bigwedge^k\tilde{\mathcal{K}}^d$, we can write 
$$\mathbf{v}=\sum_{J\subseteq [d]:\vert J\vert=k}\phi_J(\mathbf{v})\mathbf{e}_J.$$
where $\phi_J(\mathbf{v})$ is the coefficient of $\mathbf{e}_J$ in the above decomposition. We define the norm on $\bigwedge^k\tilde{\mathcal{K}}^d$ by 
$$\Vert \mathbf{v}\Vert_{k-\operatorname{vec}}:=\max_{\vert J\vert=k}\vert \phi_J(\mathbf{v})\vert.$$
and we define the support of $\mathbf{v}$ by
$$\operatorname{supp}(\mathbf{v})=\{J:\phi_J(\mathbf{v})\neq 0\}.$$
We define $\bigwedge\tilde{\mathcal{K}}^d=\bigcup_{k=0}^d\bigwedge^k\tilde{\mathcal{K}}^d$. Note that this union is disjoint, and it is often called the Grassmanian. 
\end{definition}
\begin{definition}
    A $k$ dimensional measured space is a vector space $M\subseteq \tilde{\mathcal{K}}^d$ equipped with a non-zero wedge product $\det(M)\in \bigwedge^kM$. We often denote a measured space and its determninant by $(M,\det(M))$. We define $\Vert M\Vert_{ms}=\Vert \det(M)\Vert_{k-\operatorname{vec}}$.

    Let $\Delta\subseteq \tilde{\mathcal{K}}^d$ be a discrete subgroup defined by $\Delta=\operatorname{span}_{\mathcal{R}}\big\{\mathbf{v}^{(1)},\dots, \mathbf{v}^{(k)}\big\}$. Then, we define the measured space corresponding to $\Delta$ by
    $$M(\Delta)=\operatorname{span}_{\tilde{\mathcal{K}}}(\Delta).$$
    We define the determinant of $M(\Delta)$ as 
    $$\det(M(\Delta))=\mathbf{v}^{(1)}\wedge \dots\wedge \mathbf{v}^{(k)}.$$
    If $\Delta=\{0\}$, then, we define $\det\left(\{0\}\right)=1$, and $M\left(\{0\}\right)=\left(\{0\},1\right)$. Notice that the determinant is defined up to $\operatorname{SL}_d(\mathcal{O})$ and $\Vert \det(M(\Delta))\Vert$ is well defined.
\end{definition}
\begin{definition}
    A measured flag of length $k$ is a sequence of measured spaces 
    $$F=\bigg\{0=V^{(0)}<V^{(1)}<\dots<V^{(k)}=\tilde{\mathcal{K}}^d\bigg\}.$$
    We denote the space of measured flags by $\mathcal{F}_d$ and we define 
    $$\big\Vert F\Vert_{\operatorname{flag}}=\max_{0<j<k}\Vert V^{(j)}\big\Vert_{ms}.$$
    We define the trivial flag as $\Big\{\left(\{0\}, 1\right)<\left(\tilde{\mathcal{K}}^d, \mathbf{e}^{(1)}\wedge \dots \wedge \mathbf{e}^{(d)}\right)\Big\}$.
\end{definition}
\begin{remark}
    The norms $\Vert \cdot \Vert_{k-\operatorname{vec}}$, $\Vert \cdot \Vert_{ms}$ and $\Vert \cdot \Vert_{\operatorname{flag}}$ all depend on the norm chosen on $\tilde{\mathcal{K}}^d$. To prove Theorem \ref{AxcapWR}, we shall vary the norms by acting on them with $A$ instead of acting on lattices with $A$. In general, when we use some norm which is not $\Vert \cdot \Vert_{\infty}$, we shall denote the norm of a measured space by $\Vert M\Vert_{\Vert \cdot\Vert ,ms}$, the successive minima with respect to the corresponding norm by $\lambda_{j,\Vert \cdot \Vert}$, and the flag norm by $\Vert F\Vert_{\Vert \cdot \Vert,\operatorname{flag}}$. When $\Vert \cdot \Vert=\Vert \cdot \Vert_{\infty}$, we shall remove the norm subscript completely and write $\Vert M\Vert_{ms}$, $\lambda_1\dots \lambda_d$, and $\Vert F\Vert_{\operatorname{flag}}$. 
\end{remark}
We define 
$$\mathbb{R}_0^d=\Bigg\{\mathbf{v}=\sum_{i=1}^dv_i\mathbf{e}^{(i)}\in \mathbb{R}^d:\sum_{i=1}^dv_i=0\Bigg\},$$
and we define the norm on $\mathbb{R}_0^d$ by $\Vert \mathbf{a}\Vert_{\infty}=\max_i\vert a_i\vert$. We will be interested in functions $F:\mathbb{R}_0^d\rightarrow \mathcal{F}_d$ with certain properties. This will be similar to the flag valued functions appearing in \cite[Section 2]{Sol}. 
\begin{definition}
    A flag valued function $F:\mathbb{R}_0^d\rightarrow \mathcal{F}_d$ is bounded if
    $$\sup_{\mathbf{a}\in \mathbb{R}_0^d}\Vert F(\mathbf{a})\Vert_{\operatorname{flag}}<\infty.$$
\end{definition}

\subsubsection{The Minkowski Flag}
Given a norm $\Vert \cdot\Vert$ on $\tilde{\mathcal{K}}^d$, we define the ball of radius $r$ with respect to $\Vert \cdot \Vert$ by
$$B_{\Vert \cdot \Vert}(0,r)=\{\mathbf{v}\in \tilde{\mathcal{K}}^d:\Vert \mathbf{v}\Vert\leq r\}.$$
We note that the norm $\Vert \cdot\Vert$ should be non-trivial. We define the Minkowski flag of a lattice $\Lambda\in \mathcal{L}_d$ with respect to the norm $\Vert \cdot\Vert$ as
$$f_{Mink,\Vert \cdot\Vert}(\Lambda)=\Big\{M\left(\operatorname{span}(\Lambda\cap B_{\Vert \cdot\Vert}(0,r))\cap \Lambda\right):r>0\Big\}.$$
Then $\Lambda\in \mathcal{WR}_d$ if and only if $f_{Mink,\Vert \cdot\Vert_{\infty}}(\Lambda)$ is the trivial flag. For simplicity, we denote $f_{Mink}(\Lambda)=f_{Mink,\Vert \cdot \Vert_{\infty}}(\Lambda)$. 
For $\mathbf{a}\in \mathbb{R}_0^d$ we define a norm on $\tilde{\mathcal{K}}^d$ by 
$$\Vert \mathbf{v}\Vert_{\mathbf{a}}=q^{\max_{i=1,\dots, d}(\mathbf{a}_i+\log\vert v_i\vert)}.$$
Notice that when $\mathbf{a}\in \mathbb{Z}_0^d$, then this coincides with $\Vert x^{\mathbf{a}}\mathbf{v}\Vert$, where $x^{\mathbf{a}}=\operatorname{diag}(x^{\mathbf{a}_1},\dots, x^{\mathbf{a}_d})\in A$. Define the function $F:\mathbb{R}_0^d\rightarrow \mathcal{F}_d$ by
$$F(\mathbf{a})=f_{Mink,\Vert \cdot\Vert_{\mathbf{a}}}(\Lambda).$$
Then, $A\Lambda\cap \mathcal{WR}_d\neq \emptyset$ if and only if there exists $\mathbf{a}\in \mathbb{Z}_0^d$ such that $F(\mathbf{a})=f_{Mink,\Vert \cdot\Vert_{\mathbf{a}}}(\Lambda)$ is trivial. One of our main theorems of this section is the following function field analogue of \cite[Theorem 2.3]{Sol} for the special case of the Minkowski flag. 
\begin{theorem}
\label{TrivFlag}
    For every $\Lambda\in \mathcal{L}_d$, there exists $\mathbf{a}\in \mathbb{R}_0^d$ such that $F(\mathbf{a})$ is trivial.
\end{theorem}
Thus, in view of Theorem \ref{TrivFlag} to prove Theorem \ref{AxcapWR}, it suffices to show that we can choose $\mathbf{a}\in \mathbb{Z}_0^d$ such that $F(\mathbf{a})$ is trivial. 
\begin{theorem}
\label{ainZ_0^d}
    There exists $\mathbf{a}\in \mathbb{Z}_0^d$ such that $F(\mathbf{a})$ is trivial. As a consequence, $A\Lambda\cap \mathcal{WR}_d\neq \emptyset$.
\end{theorem}
To prove Theorem \ref{ainZ_0^d}, we shall prove the following proposition, which establishes a connection between the norms $\Vert \cdot \Vert_{\mathbf{a}}$ and convex bodies.
\begin{proposition}
\label{F(a)TrivCond}
    Let $\mathbf{a}\in \mathbb{R}_0^d$. 
    If $F(\mathbf{a})$ is trivial, then $\mathbf{a}_1,\dots, \mathbf{a}_d$ are equal mod $1$, where $\mathbf{a}=(\mathbf{a}_1\dots \mathbf{a}_d)$. This implies that there exists a convex body $\mathcal{C}$ and a constant $c>0$, such that $\|\cdot\|_{\mathcal{C}} = c\|\cdot\|_{\mathbf{a}}$. 
\end{proposition}
\begin{proof}
    Let $\mathbf{a}\in \mathbb{R}_0^d$ satisfy that $F(\mathbf{a})$ is trivial, that is, $\Lambda$ satisfies that \[\lambda_{1,\|\cdot\|_{\mathbf{a}}}(\Lambda) = \lambda_{2,\|\cdot\|_{\mathbf{a}}}(\Lambda) = \dots = \lambda_{d,\|\cdot\|_{\mathbf{a}}}(\Lambda).\]
    Assume on the contrary that $\mathbf{a}_1,\dots, \mathbf{a}_d$ are not all equal mod $1$.
    Denote $r=\lambda_{1,\Vert \cdot \Vert_{\mathbf{a}}}(\Lambda)$.
    Since $F(\mathbf{a})$ is trivial, then we deduce that $B_{\|\cdot\|_\mathbf{a}}(0, r)\cap \Lambda$ contains $d$ linear independent vectors $\mathbf{v}^{(1)}\dots \mathbf{v}^{(d)}$, while for every $r'<r$ we have 
    \begin{align}\label{eq: ball empty}
        B_{\|\cdot\|_\mathbf{a}}(0, r')\cap\Lambda = \{0\}. 
    \end{align}
    The image of $\|\cdot\|_{\mathbf{a}}:\tilde{\mathcal{K}}^d\rightarrow \mathbb{R}$ is 
    \[{\rm  Im}(\|\cdot\|_{\mathbf{a}}) = \{0\} \cup \{q^s:s\in \mathbf{Z} + \mathbf{a}_i,i=1\dots d\}.\]
    Suppose that $r = q^{\ba_i + n}$ for some $n\in \mathbf{Z}$ and some $1\le i\leq d$. Let $r' = q^{\mathbf{a}_j+m}<r$ be the maximal value in ${\rm  Im}(\|\cdot\|_{\ba})$ that is smaller then $r$. Since $\mathbf{a}_1,\dots, \mathbf{a}_d$ are not all equal mod $1$, then by choosing $m$ properly, we can ensure that $\mathbf{a}_i+n-(\mathbf{a}_j+m)>-1$ whenever $\mathbf{a}_j\neq \mathbf{a}_i\mod 1$. Hence, $r' > \frac{r}{q}$. Since there exists some $\mathbf{v}\in \tilde{\mathcal{K}}^d$ such that $\Vert \mathbf{v}\Vert_{\mathbf{a}}=r'$, we deduce that
    \[B_{\|\cdot\|_\mathbf{a}}\left(0, \frac{r}{q}\right)\subsetneq B_{\|\cdot\|_\mathbf{a}}(0, r')\subsetneq B_{\|\cdot\|_\mathbf{a}}(0, r).\]
    Let $\sum_{i=1}^d\alpha_i\mathbf{v}^{(i)}\in B_{\Vert \cdot \Vert_{\mathbf{a}}}(0,r)$. By writing $\alpha_i=[\alpha_i]+\langle \alpha_i\rangle$, where $\vert\langle \alpha_i\rangle\vert\leq q^{-1}$ and $[\alpha_i]\geq 1$, we obtain that $\sum_{i=1}^d\langle\alpha_i\rangle\mathbf{v}^{(i)}\in B_{\Vert \cdot \Vert_{\mathbf{a}}}\left(0,\frac{r}{q}\right)$. Hence,
    $$B_{\|\cdot\|_\mathbf{a}}(0, r) / B_{\|\cdot\|_\mathbf{a}}\left(0, \frac{r}{q}\right) \cong \mathcal{R}^d,$$
    by identifying a coset $\sum_{i=1}^d\alpha_i\mathbf{v}^{(i)}+B_{\Vert \cdot \Vert_{\mathbf{a}}}\left(0,\frac{r}{q}\right)$ with the vector $([\alpha_1],\dots,[\alpha_d])$. Thus, 
    $$\dim_{\mathcal{R}}\left(B_{\|\cdot\|_\mathbf{a}}(0, r) / B_{\|\cdot\|_\mathbf{a}}(0, r')\right)<d.$$
    Hence the images of $\mathbf{v}^{(1)}\dots \mathbf{v}^{(d)}$ in 
    $$B_{\|\cdot\|_\mathbf{a}}(0, r) / B_{\|\cdot\|_\mathbf{a}}(0, r')$$
    are linearly dependent. Thus, there exist $b_1,\dots,b_n\in \mathcal{R}$ such that 
    $$b_1\mathbf{v}^{(1)}+b_2\mathbf{v}^{(2)}+\dots +b_d\mathbf{v}^{(2)} \in B_{\|\cdot\|_\mathbf{a}}(0, r').$$
    Then, $b_1\mathbf{v}^{(1)}+b_2\mathbf{v}^{(2)}+\dots +b_d\mathbf{v}^{(d)} \in B_{\|\cdot\|_\mathbf{a}}(0, r') \cap \Lambda$. Since $\mathbf{v}^{(1)}\dots \mathbf{v}^{(d)}$ are linearly independent, then $b_1\mathbf{v}^{(1)}+\dots +b_d\mathbf{v}^{(d)}\neq 0$. This contradicts (\ref{eq: ball empty}), which implies that $\ba_1,\dots, \ba_n$ are equal mod $1$. Let $c=q^{\mathbf{a}_1\mod 1}$ and let $\mathcal{C} = B_{\|\cdot\|_\mathbf{a}}(0, c)$. Then $\operatorname{Im}(\Vert \cdot\Vert_{\mathbf{a}})=\{0\}\cup cq^{\mathbb{Z}}$ and $\|\cdot\|_{\mathbf{a}}=c\|\cdot\|_{\mathcal{C}}$.
\end{proof}
We shall now use Theorem \ref{TrivFlag} and Proposition \ref{F(a)TrivCond} to prove Theorem \ref{ainZ_0^d}.
\begin{proof}[Proof of Theorem \ref{ainZ_0^d}]
Assume that $\Lambda\notin \mathcal{WR}_d$. Otherwise, there is nothing to prove. Let $\mathbf{a}=\mathbf{a}^{(1)}+\mathbf{a}^{(2)}$ satisfy the conclusion of Theorem \ref{TrivFlag}, where $\mathbf{a}^{(1)}\in \mathbb{Z}_0^d$ and $\max_i\mathbf{a}^{(2)}_i<1$. First assume that $\mathbf{a}^{(1)}=0$. Let $\mathbf{v}^{(1)}\dots \mathbf{v}^{(d)}$ be a set of successive minima vectors for $\Lambda$ with respect to $\Vert \cdot \Vert_{\mathbf{a}}=\Vert \cdot\Vert_{\mathbf{a}^{(2)}}$. Hence, by Proposition \ref{F(a)TrivCond} there exist $c>0$ and a convex body $\mathcal{C}=B_{\Vert \cdot \Vert_{\mathbf{a}^{(2)}}}(0,c)$ such that $c\Vert \cdot \Vert_{\mathcal{C}}=\Vert \cdot \Vert_{\mathbf{a}^{(2)}}$. Then, $\operatorname{Vol}\left(B_{\Vert \cdot \Vert_{\mathbf{a}}^{(2)}}\left(0,c\right)\right)=\operatorname{Vol}(\mathcal{C})=c^d$ and for every $i=1\dots d$, $\lambda_{i,\Vert \cdot\Vert_{\mathbf{a}^{(2)}}}(\Lambda)=c\lambda_{i,\mathcal{C}}(\Lambda)$. Thus by Theorem \ref{Mink2ndConv},
\begin{equation}
    1=c^d\prod_{i=1}^d\lambda_{i,\mathcal{C}}(\Lambda)=\prod_{i=1}^d\lambda_{i,\Vert \cdot \Vert_{\mathbf{a}^{(2)}}}(\Lambda).
\end{equation}
Thus, if $F(\mathbf{a}^{(2)})$ were trivial, then, $\lambda_{1,\Vert \cdot\Vert_{\mathbf{a}^{(2)}}}(\Lambda)=1$. On the other hand, since $\Lambda\notin \mathcal{WR}_d$, then there exists $\mathbf{v}\in \Lambda$ such that $\Vert \mathbf{v}\Vert\leq q^{-1}$. Hence, for every $i=1\dots d$, $\mathbf{a}_i^{(2)}+\log\vert v_i\vert<1-1=0$, and therefore, 
$$\big\Vert \mathbf{v}\big\Vert_{\mathbf{a}^{(2)}}=q^{\max_{i=1\dots d}\left(\mathbf{a}_i+\log\vert v_i\vert\right)}<1.$$
Thus, $F(\mathbf{a}^{(2)})$ cannot be trivial. Therefore, $\mathbf{a}^{(1)}\neq 0$. 

On the other hand, if $F\left(\mathbf{a}^{(1)}\right)$ is not trivial, then by the above claim, by replacing $F$ with the function 
$$\tilde{F}(\mathbf{a}')=F_{Mink,\Vert \cdot\Vert_{\mathbf{a}'}}\left(x^{\mathbf{a}^{(1)}}\Lambda\right)=F_{Mink,\Vert \cdot\Vert_{\mathbf{a}'+\mathbf{a}^{(1)}}}(\Lambda),$$
we obtain that $\tilde{F}\left(\mathbf{a}^{(2)}\right)=F\left(\mathbf{a}^{(1)}+\mathbf{a}^{(2)}\right)=F(\mathbf{a})$ is also not trivial which is a contradiction to the assumption that $F(\mathbf{a})$ is trivial. Hence $F\left(\mathbf{a}^{(1)}\right)$ must be trivial. As a consequence, $x^{\mathbf{a}^{(1)}}\Lambda\in \mathcal{WR}_d$, and hence, $A\Lambda\cap \mathcal{WR}_d\neq \emptyset$. 
\end{proof}
This concludes the proof of Theorem \ref{AxcapWR}. Thus, it suffices to prove Theorem \ref{TrivFlag}. We shall now show that the Minkowski flag satisfies certain properties. Firstly, we will need an analogue of \cite[Lemma 2.4]{Sol}. This is a reinterpretation of the upper bound in Theorem \ref{Mink2nd} in the language of flags.
\begin{lemma}
\label{fMinkBndd}
For any $\Lambda\in \mathcal{L}_d$, $\big\Vert f_{Mink}(\Lambda)\big\Vert_{\operatorname{flag}}\leq 1$. 
\end{lemma}
\begin{proof}
    By Theorem \ref{Mink2nd}, for every $\Lambda\in \mathcal{L}_d$, $\prod_{i=1}^d\lambda_i(\Lambda)=1$. Hence, for every $k=1\dots d$, $\prod_{i=1}^k\lambda_i(\Lambda)\leq 1$. Let $\mathbf{v}^{(1)}\dots \mathbf{v}^{(d)}\in \Lambda$ be a basis of successive minima for $\Lambda$. Let $V=M\left(\operatorname{span}(\Lambda\cap B(0,r))\cap \Lambda\right)$ be of dimension $k$. We shall show that $\Vert V\Vert_{ms}\leq 1$. Then, by definition of the successive minima $\lambda_k(\Lambda)\leq r<\lambda_{k+1}(\Lambda)$. Hence, $V=\operatorname{span}_{\tilde{\mathcal{K}}}\big\{\mathbf{v}^{(1)}\dots \mathbf{v}^{(k)}\}$ and thus, $\det(V)=\mathbf{v}^{(1)}\wedge \dots \wedge \mathbf{v}^{(k)}$. Notice that if $\big\Vert\mathbf{u}^{(j)}\big\Vert \leq 1$ for every $j=1\dots d$, then $\big\Vert \mathbf{u}^{(1)}\wedge \dots \wedge \mathbf{u}^{(k)}\big\Vert\leq 1$. Hence,
    \begin{equation}
        \big\Vert \mathbf{v}^{(1)}\wedge\dots \wedge\mathbf{v}^{(k)}\big\Vert\leq \lambda_1(\Lambda)\dots \lambda_k(\Lambda)\Bigg\Vert \frac{\mathbf{v}^{(1)}}{\lambda_1(\Lambda)}\wedge \dots \wedge \frac{\mathbf{v}^{(k)}}{\lambda_k(\Lambda)}\Bigg\Vert\leq 1.
        \end{equation}
Hence, $\Vert V\Vert_{ms}=\big\Vert \mathbf{v}^{(1)}\wedge \dots \wedge \mathbf{v}^{(k)}\Vert\leq 1$, and therefore, $\big\Vert f_{Mink}(\Lambda)\big\Vert_{\operatorname{flag}}\leq 1$.
\end{proof}
\subsubsection{Proof of Theorem \ref{TrivFlag}}
In order to prove Theorem \ref{TrivFlag}, we shall use \cite[Theorem 1.4]{Sol}. We shall first define the invariance dimension of a convex open set $U\subseteq \mathbb{R}^d$.
\begin{definition}
    It is well known that $\mathbb{R}^d$ acts on the set of convex open sets in $\mathbb{R}^d$ by translation. The invariance dimension of a convex open set $U\subseteq \mathbb{R}^d$ is the dimension of its stabilizer under this action. We denote
    $$\operatorname{invdim}(U)=\dim\operatorname{stab}_{\mathbb{R}^d}(U),$$
    and we define $\operatorname{invdim}(\emptyset)=-\infty$. Given a set $U\subseteq \mathbb{R}^d$, we denote its convex hull by $\operatorname{conv}(U)$. 
\end{definition} 
\begin{theorem}{\cite[Theorem 1.4]{Sol}}
    \label{covClaim}
    Let $\mathcal{U}$ be an open cover of $\mathbb{R}^d$. Assume that 
    \begin{enumerate}
        \item \label{LocFin} The cover $\{\operatorname{conv}(U):U\in \mathcal{U}\}$ is locally finite, that is every compact set intersects at most finitely many elements of the cover $\{\operatorname{conv}(U):U\in \mathcal{U}\}$.
        \item \label{BnddInvDim} For every $k\leq d$ and for every $U_1\dots U_k\in \mathcal{U}$, 
        $$\operatorname{invdim conv}(U_1\cap \dots \cap U_k)\leq d-k.$$
    \end{enumerate}
    Then, there are $d+1$ sets in $\mathcal{U}$ with non-trivial intersection.
\end{theorem}
We shall now prove Theorem \ref{TrivFlag}. Our proof is absolutely analogous to the proof of Theorem 2.3 in \cite{Sol}, but we shall include it to portray the differences and similarities between the real and function field settings. Unlike the real setting, in which Solan used the action of the diagonal group on the Grassmanian, here we look at the $A$ action on the norm on $\mathbb{R}_0^d$.  
\begin{proof}[Proof of Theorem \ref{TrivFlag}]
    Assume to the contrary that $F(\mathbf{a})$ is nowhere trivial, that is for every $\mathbf{a}\in \mathbb{R}_0^d$, there exists some $\mathbf{v}\in F(\mathbf{a})$ such that $0<\mathbf{v}<\tilde{\mathcal{K}}^d$. For every $\mathbf{v}\in \bigwedge\tilde{\mathcal{K}}^d$ with $1\leq \dim (\mathbf{v})\leq d-1$, define
    $$U_{\mathbf{v}}=\{\mathbf{a}\in \mathbb{R}_0^d:\mathbf{v}\in F(\mathbf{a})\}.$$
    We first note that due to the definition of $F$, $U_{\mathbf{v}}$ is non-empty if and only if $\mathbf{v}\in \bigwedge\Lambda$. Since $F(\mathbf{a})$ is nowhere trivial, then $\mathcal{U}=\{U_{\mathbf{v}}\}_{\mathbf{v}\in \bigwedge\tilde{\mathcal{K}}^d,1\leq \dim(\mathbf{v})\leq d-1}=\{U_{\mathbf{v}}\}_{\mathbf{v}\in \bigwedge\Lambda,1\leq \dim(\mathbf{v})\leq d-1}$ is a cover of $\mathbb{R}_0^d\cong \mathbb{R}^{d-1}$, and we shall now show that the cover $\mathcal{U}$ satisfies the conditions of Theorem \ref{covClaim}. 

    We shall first show that $\mathcal{U}$ is an open cover. Firstly, notice that if $\mathbf{a}\in U_{\mathbf{v}}$, then, there exists some $r>0$ such that $\mathbf{v}=M\left(\operatorname{span}(\Lambda\cap B_{\Vert \cdot \Vert_{\mathbf{a}}}(0,q^r)\right)\cap \Lambda)$. By definition 
    $$B_{\Vert \cdot \Vert_{\mathbf{a}}}(0,q^r)=\Big\{\mathbf{u}\in \tilde{\mathcal{K}}^d:\max_iq^{a_i+\log\vert u_i\vert}<q^r\Big\}.$$
    Since $B_{\Vert \cdot \Vert_{\mathbf{a}}}(0,q^r)\cap \Lambda$ is a finite set, then by varying $r\in \mathbb{R}$ while preserving $B_{\Vert \cdot\Vert_{\mathbf{a}}}(0,q^r)\cap\Lambda$, we can find some $r_1<r_2$ such that 
    $$B_{\Vert \cdot\Vert_{\mathbf{a}}}(0,q^{r_1})\cap\Lambda=B_{\Vert \cdot\Vert_{\mathbf{a}}}(0,q^{r_2})\cap\Lambda,$$
    and 
    $$\mathbf{v}=M\left(\operatorname{span}\left(B_{\Vert \cdot\Vert_{\mathbf{a}}}(0,q^{r_1}\right)\cap\Lambda)\cap \Lambda\right)=M\left(\operatorname{span}(B_{\Vert \cdot \Vert_{\mathbf{a}}}(0,q^{r_2})\cap \Lambda)\cap \Lambda\right).$$
    Hence, if $\Vert \mathbf{a}'-\mathbf{a}\Vert<r_2-r_1$, then, 
    \begin{equation}
        \begin{split}
            B_{\Vert \cdot \Vert_{\mathbf{a}'}}(0,q^{r_1})=\Big\{\mathbf{u}\in \tilde{\mathcal{K}}^d:\max_iq^{a_i'+\log\vert u_i\vert}<q^{r_1}\Big\}\\
        \subseteq \Big\{\mathbf{u}\in \tilde{\mathcal{K}}^d:\max_iq^{a_i+\log \vert u_i\vert}<q^{r_2}\Big\}=B_{\Vert \cdot \Vert_{\mathbf{a}}}(0,q^{r_2}).
        \end{split}
    \end{equation}
    Similarly, by replacing the roles of $\mathbf{a}$ and $\mathbf{a}'$,
    \begin{equation}
        \begin{split}
            B_{\Vert \cdot\Vert_{\mathbf{a}}}(0,q^{r_1})=\Big\{\mathbf{u}\in \tilde{\mathcal{K}}^d:\max_iq^{a_i+\log\vert u_i\vert}<q^{r_1}\Big\}\\
            \subseteq \Big\{\mathbf{u}\in \tilde{\mathcal{K}}^d:\max_iq^{a_i'+\log\vert u_i\vert}<q^{r_2}\Big\}=B_{\Vert \cdot \Vert_{\mathbf{a}'}}(0,q^{r_2}).
        \end{split}
    \end{equation}
    Hence, 
    $$B_{\Vert \cdot \Vert_{\mathbf{a}'}}(0,q^{r_1})\cap \Lambda\subseteq B_{\Vert \cdot \Vert_{\mathbf{a}}}(0,q^{r_2})\cap \Lambda,$$
    and 
    $$B_{\Vert \cdot \Vert_{\mathbf{a}}}(0,q^{r_1})\cap \Lambda\subseteq B_{\Vert \cdot \Vert_{\mathbf{a}'}}(0,q^{r_2})\cap \Lambda.$$
    Since $B_{\Vert \cdot \Vert_{\mathbf{a}'}}(0,q^r)\cap \Lambda$ is finite then there exists some $r_1\leq r'\leq r_2$ such that 
    $$B_{\Vert \cdot \Vert_{\mathbf{a}'}}(0,q^{r'})\cap \Lambda=B_{\Vert \cdot \Vert_{\mathbf{a}}}(0,q^{r_1})\cap \Lambda.$$
    Thus, 
    $$\mathbf{v}=M\left(\operatorname{span}\left(B_{\Vert \cdot \Vert_{\mathbf{a}}}(0,q^{r_1}\cap \Lambda\right)\cap \Lambda\right)=M\left(\operatorname{span}\left(B_{\Vert \cdot \Vert_{\mathbf{a}'}}(0,q^{r'}\cap \Lambda\right)\cap \Lambda\right).$$
    Thus, $\mathbf{a}'\in U_{\mathbf{v}}$, so that $B(\mathbf{a},r_2-r_1)\subseteq U_{\mathbf{v}}$. Thus, $\mathcal{U}$ is an open cover. 
    
    We shall now prove that $\mathcal{U}$ satisfies condition (\ref{LocFin}). Let  
    $$U_{\mathbf{v}}'=\Big\{\mathbf{a}\in \mathbb{R}_0^d:\Vert \mathbf{v}\Vert_{\Vert \cdot \Vert_{\mathbf{a}},{ms}}\leq 1\Big\},$$
    Since $F$ is bounded by $1$, then, $U_{\mathbf{v}}\subseteq U_{\mathbf{v}}'$ for every $\mathbf{v}$ such that $U_{\mathbf{v}}\neq \emptyset$. Let $\mathcal{U}'=\{U_{\mathbf{v}}':\mathbf{v}\in \bigwedge\Lambda\}$.  We first prove that $\mathcal{U}'$ satisfies condition (\ref{LocFin}), which will imply that $\mathcal{U}$ satisfies condition (\ref{LocFin}). Since $\Lambda$ is discrete, then $\bigwedge\Lambda$ is also discrete. Thus, $\bigwedge \Lambda$ is locally finite, that is its intersection with every compact set in $\bigwedge\tilde{\mathcal{K}}^d$ is finite. Notice that if $X\subseteq \mathbb{R}_0^d$ is compact, then, there exists some $r>0$, such that $X\subseteq \{\mathbf{a}\in \mathbb{R}_0^d:\max_{i=1,\dots,d} a_i\leq r\}$. Thus, for every $\mathbf{a}\in X$ and for every $\mathbf{v}\in \bigwedge\tilde{\mathcal{K}}^d$, we have 
    \begin{equation}
    \label{eqn:aContract}
        \Vert \mathbf{v}\Vert_{\Vert \cdot \Vert_{\mathbf{a}},ms}=q^{\max_{i=1,\dots,d}(a_i+\log_q\vert v_i\vert)}\geq q^{-(d-1)r}\Vert \mathbf{v}\Vert_{ms}.
    \end{equation}
    Thus, if $\Vert\mathbf{v}\Vert_{\Vert\cdot\Vert_{\mathbf{a}},ms}\leq 1$, then, $\Vert \mathbf{v}\Vert_{ms}\leq q^{(d-1)r}$. Since $\bigwedge\Lambda$ is discrete, then, there are finitely many $\mathbf{v}\in \bigwedge\Lambda$ satisfying $\Vert \mathbf{v}\Vert\leq q^{(d-1)r}$. Hence, $X$ intersects finitely many sets in $\mathcal{U}'$, so that $\mathcal{U}'$ is finite. 

    We shall now show that $\mathcal{U}$ satisfies condition (\ref{BnddInvDim}). Let $U_{\mathbf{v}^{(1)}}\dots U_{\mathbf{v}^{(l)}}\in \mathcal{U}$ be distinct elements with 
    $$U_{\mathbf{v}^{(1)}}\cap \dots \cap U_{\mathbf{v}^{(l)}}=V\neq \emptyset.$$
    Then, for every $\mathbf{a}\in V$, we have $\mathbf{v}^{(1)},\dots, \mathbf{v}^{(l)}\in F(\mathbf{a})$. Hence, $\mathbf{v}^{(1)},\dots, \mathbf{v}^{(l)}$ form a flag. Assume without loss of generality that $$0<\mathbf{v}^{(1)}<\dots <\mathbf{v}^{(l)}<\tilde{\mathcal{K}}^d.$$
    By \cite[Lemma 2.5]{Sol}, there exists a permutation $\sigma:[d]\rightarrow [d]$, such that $\sigma\left(\left[\dim \mathbf{v}^{(j)}\right]\right)\in \operatorname{supp}\left(\mathbf{v}^{(j)}\right)$ for every $j=1\dots l$, where $[m]=\{1\dots m\}$. Assume without loss of generality that $\sigma=Id$. Denote $c_j:=\left|\varphi_{\left[\dim \left(\mathbf{v}^{(j)}\right)\right]}\left(\mathbf{v}^{(j)}\right)\right|$. Then, for any $\mathbf{a}\in V$, 
    \begin{equation}
        1>\Vert F(\mathbf{a})\Vert_{\operatorname{flag}}\geq \max_{j=1,\dots l} \big\Vert \mathbf{v}^{(j)}\big\Vert_{\Vert\cdot\Vert_{\mathbf{a}},ms}\geq \max_{j=1\dots l}c_jq^{\sum_{i\in \left[\dim\left(\mathbf{v}^{(j)}\right)\right]}a_i}.
    \end{equation}
    Hence, $V$ is contained in the convex set 
    \begin{equation}
        P=\bigcap_{j=1}^l\psi_{\dim\left(\mathbf{v}^{(j)}\right)}^{-1}(-\infty,-\log_q c_j),
    \end{equation}
    where $\psi_m:\mathbb{R}_0^d\rightarrow \mathbb{R}$ is defined by $\psi_m(\mathbf{a})=\sum_{i=1}^ma_i$. Hence, $\operatorname{conv}(V)\subseteq P$. Since $\psi_{\dim\left(\mathbf{v}^{(j)}\right)}$ are independent for $j=1\dots l$, then, $\operatorname{invdim}(P)=d-1-l$. Hence, $\operatorname{invdim}(V)\leq d-1-l$, which proves that condition (\ref{BnddInvDim}) holds. 

    Thus, by Theorem \ref{covClaim}, there exist $d$ sets $U_{\mathbf{v}^{(1)}}\dots U_{\mathbf{v}^{(d)}}$ with nontrivial intersection. Thus, there exists $\mathbf{a}\in U_{\mathbf{v}^{(1)}}\cap \dots\cap U_{\mathbf{v}^{(d)}}$. Therefore, $\mathbf{v}^{(1)},\dots,\mathbf{v}^{(d)}\in F(\mathbf{a})$. Hence, they form a flag of the form 
    $$\{0\}<\mathbf{v}^{(1)}<\mathbf{v}^{(2)}<\dots<\mathbf{v}^{(d)}<\tilde{\mathcal{K}}^d,$$
    and this is impossible. Hence, there exists some $\mathbf{a}\in \mathbb{R}_0^d$ such that $F(\mathbf{a})$ is trivial. 
\end{proof}
\subsection{Proof of Theorem \ref{Mink2}}
We shall now deduce Theorem \ref{Mink2} from Theorem \ref{MinkConj} and the following theorem by Aggarwal (see \cite{Agg}). 
\begin{theorem}
    \label{Agg}
    For any $\theta\in \tilde{\mathcal{K}}$, $\mu\left(\begin{pmatrix}
        1&\theta\\
        0&1
    \end{pmatrix}\mathcal{R}^2\right)=q^{-2}$.
\end{theorem}
To conclude the proof of Theorem \ref{Mink2}, We shall show that almost every lattice of a certain form has a dense $A$ orbit. This is a well known theorem which we shall prove for completeness. For the real analogue see \cite[Lemma 4.8]{S09}. 
\begin{lemma}
\label{DenseUni}
Let $d\geq 2$, and let $I_{d-1}$ denote the identity matrix of order $d-1\times d-1$. Then for almost each $\boldsymbol{\theta}\in \tilde{\mathcal{K}}^{d-1}$, the lattice 
\begin{equation}
\label{eqn:uniLatt}
    \Lambda=\begin{pmatrix}
I_{d-1}&\boldsymbol{\theta}\\
0&1
\end{pmatrix}\mathcal{R}^d
\end{equation}
has a dense $A$ orbit. 
\end{lemma}
\begin{proof}[Proof of Theorem \ref{Mink2}]
    By Lemma \ref{DenseUni}, for almost every $\theta\in \tilde{\mathcal{K}}$, $\Lambda_{\theta}$ has a dense $A$ orbit. Hence, by upper semicontinuity of $\mu$, ergodicity of the $A$ action \cite{Mau}, and Theorem \ref{Agg}, for almost every $\theta \in \tilde{\mathcal{K}}$, $q^{-2}=\mu(\Lambda_{\theta})=\min\mathcal{S}_2$. On the other hand, by Theorem \ref{MinkConj}, for every $\Lambda\in \mathcal{L}_2$, $\mu(\Lambda)\leq q^{-2}$. Hence $\mathcal{S}_2=\{q^{-2}\}$.
    \end{proof}
\begin{proof}[Proof of Lemma \ref{DenseUni}]
For $\boldsymbol{\theta}\in \tilde{\mathcal{K}}^{d-1}$ let $u_{\boldsymbol{\theta}}=\begin{pmatrix}
1&\boldsymbol{\theta}\\
0&1
\end{pmatrix}$ and let $\Lambda_{\boldsymbol{\theta}}=u_{\boldsymbol{\theta}}\mathcal{R}^d$. For $k\geq 0$, let 
$$a_k=\operatorname{diag}\big\{x^k,\dots x^k,x^{-nk}\big\}$$
It is well known that 
$$U^+(a_k)=\big\{g\in G:\lim_{k\rightarrow \infty}a_{-k}ga_k=I\big\}=\big\{u_{\boldsymbol{\theta}}:\boldsymbol{\theta}\in \tilde{\mathcal{K}}^{d-1}\big\}.$$
Denote
$$U^-(a_k)=\big\{g\in G:\lim_{k\rightarrow \infty}a_kga_{-k}=I\big\},$$
$$U^0(a_k)=\{g\in G:\forall k\in \mathbb{Z}, a_kg=ga_k\}.$$
Then, for each $\Lambda\in \mathcal{L}_d$, there exist some neighborhoods of the corresponding identity elements, $W_{\Lambda}^+\subseteq U^+(a_k)$, $W_{\Lambda}^-\subseteq U^-(a_k)$ and $W_{\Lambda}^0\subseteq U^0(a_k)$ such that the map $W_{\Lambda}^+\times W_{\Lambda}^-\times W_{\Lambda}^0\mapsto \mathcal{L}_d$ given by $(c,g,u_{\boldsymbol{\theta}})\mapsto cgu_{\boldsymbol{\theta}}\Lambda$ is a homeomorphism with image $W_{\Lambda}\subseteq \mathcal{L}_d$. Notice that if $\Lambda'=cgu_{\boldsymbol{\theta}}\Lambda\in W_{\Lambda}$, then, as $k\rightarrow \infty$,
Thus, $\{a_k\Lambda'\}_{k\geq 0}$ is dense in $\mathcal{L}_d$ if and only if $\{a_ku_{\boldsymbol{\theta}}\Lambda\}_{k\geq 0}$ is dense in $\mathcal{L}_d$. Since the one parameter group $\{a_k\}$ acts ergodically on $\mathcal{L}_d$, then, for almost every $\Lambda'\in W_{\Lambda}$, $\{a_k\Lambda'\}_{k\geq 0}$ is dense in $\mathcal{L}_d$. Therefore, for almost every $\boldsymbol{\theta}$ such that $u_{\boldsymbol{\theta}}\in U^+(a_k)$, $\{a_ku_{\boldsymbol{\theta}}\Lambda\}_{k\geq 0}$ is dense in $\mathcal{L}_d$. 

We now view $U^+(a_k)$ as a subset of $\tilde{\mathcal{K}}^{d-1}$ by identifying $u_{\boldsymbol{\theta}}$ with $\boldsymbol{\theta}$. To finish the claim, we take a countable union $\boldsymbol{\theta}_k\in \tilde{\mathcal{K}}^{d-1}$ such that $\tilde{\mathcal{K}}^{d-1}=\cup_{k\in \mathbb{N}}\left(\boldsymbol{\theta}_k+W_{\Lambda_{\boldsymbol{\theta}_k}}^+\right)$. Then,
$$W_{\Lambda_{\boldsymbol{\theta}_k}}^+\Lambda_{\boldsymbol{\theta}_k}=\{\Lambda_{\boldsymbol{\eta}}:\boldsymbol{\eta}\in \boldsymbol{\theta}_k+W_{\boldsymbol{\theta}_k}^+\}.$$ 
Thus, for almost every $\boldsymbol{\theta}\in \tilde{\mathcal{K}}^{d-1}$, $\overline{A\Lambda_{\boldsymbol{\theta}}}$ is dense in $\mathcal{L}_d$.
\end{proof}
\subsection{Computations of \texorpdfstring{$\mu$}{} for Several Lattices}
\label{subsec:compMu}
Since we were not able to find lattices $\Lambda\in \mathcal{L}_d$ with $\mu(\Lambda)<q^{-d}$ for $d\geq 3$,  we used computer experiments to attempt to find such lattices. We focused on lattices of the form
$$\Lambda_{\theta,\phi}=\begin{pmatrix}
    1&0&\theta\\
    0&1&\phi\\
    0&0&1
\end{pmatrix}\mathcal{R}^3,$$
where $\theta,\phi\in \mathcal{K}$, since $\mu(\Lambda_{\theta,\phi})$ can be computed in finitely many steps. We note that for such a lattice
\begin{equation}
    \label{eqn:mu(Lambda_theta,phi)}
    \mu(\Lambda_{\theta,\phi})=\sup_{\alpha,\beta,\gamma\in \tilde{\mathcal{K}}}\min\Big\{\vert \alpha\beta\gamma\vert,\inf_{0\neq N\in \mathcal{R}}\vert N\vert \cdot\vert \langle N\theta-\alpha\rangle\vert\cdot\vert\langle N\phi-\beta\rangle\vert\Big\},
\end{equation}
where $\Big\langle\sum_{n=-\infty}^Na_nx^n\Big\rangle=\sum_{n=-\infty}^{-1}a_nx^n$. Firstly, if indeed $\mu(\Lambda_{\theta,\phi})\leq q^{-4}$, then, for every $\vert \alpha\vert=\vert\beta\vert=\vert \gamma\vert=\frac{1}{q}$, 
\begin{equation}
    \label{eqn:mu<q**-4}
    \inf_{0\neq N\in \mathcal{R}}\vert N\vert \cdot\vert \langle N\theta-\alpha\rangle\vert\cdot\vert\langle N\phi-\beta\rangle\vert\leq q^{-4}<\vert \alpha\beta\gamma\vert=q^{-3}.
\end{equation}
Hence, it suffices to consider grids $\Lambda_{\theta,\phi}+\begin{pmatrix}
\alpha\\
\beta\\
\gamma\end{pmatrix}$ where $\vert \alpha\vert=\vert \beta\vert=\vert \gamma\vert=\frac{1}{q}$.

\begin{definition}
    Let $\alpha\in \mathcal{K}$. We say that $\operatorname{denom}(\alpha)=h$ if there exist $\alpha_j\in \mathbb{F}_q$ such that $\langle\alpha\rangle=\sum_{j=1}^h\alpha_jx^{-j}$ and $\alpha_h\neq 0$. Notice that if $\operatorname{denom}(\alpha)\leq h$, then $\vert \alpha\vert\geq q^{-h}$. 
\end{definition}
Moreover, if $\theta,\phi\in \mathcal{K}$ satisfy $\operatorname{denom}(\theta)\leq h$ and $\operatorname{denom}(\phi)\leq h$, then to compute (\ref{eqn:mu<q**-4}), it suffices to consider $N\in \mathcal{R}$ with $\vert N\vert\leq q^{h-1}$. Notice that under these conditions, $\operatorname{denom}(\langle N\theta\rangle)\leq h$ and $\operatorname{denom}(\langle N\phi\rangle)\leq h$. Hence, it suffices to consider $\alpha,\beta$ with $\operatorname{denom}(\alpha),\operatorname{denom}(\beta)\leq h$. 

We shall now reinterpret (\ref{eqn:mu<q**-4}) in a different form. Assume that $\theta,\phi,\alpha,\beta$ all have denominator which is at most $h$ and that $\vert N\vert=q^m<q^h$. Then, 
$$\vert N\vert \cdot \vert \langle N\theta-\alpha\rangle\vert\cdot\vert\langle N\phi-\beta\rangle\vert\leq q^{-4},$$
if and only if 
$$\vert\langle N\theta-\alpha\rangle\vert\cdot \vert \langle N\phi-\beta\rangle\vert\leq q^{-(m+4)}.$$
This happens if and only if there exists some $1\leq l\leq m+4$ such that both of the following conditions hold:
\begin{equation}
\label{eqn:EachCase}
\begin{cases}
    \vert \langle N\theta-\alpha\rangle\vert \leq q^{-l}\\
    \vert \langle N\phi-\beta\rangle\vert \leq q^{-(m-l+4)}
\end{cases}.
\end{equation}
Let $N=\sum_{j=0}^mN_jx^j$, $\theta=\sum_{j=1}^h\theta_jx^{-j}$, $\phi=\sum_{j=1}^h\phi_jx^{-j}$, $\alpha=\sum_{j=1}^h\alpha_jx^{-j}$, and $\beta=\sum_{j=1}^h\beta_jx^{-j}$. Then, if we write out $\langle N\theta\rangle=\sum_{i=1}^h(N\theta)_ix^{-i}$, then,
\begin{equation}
    (N\theta)_i=\sum_{j=0}^{m}N_j\theta_{i+j}.
\end{equation}
Hence, (\ref{eqn:EachCase}) holds if and only if the following conditions hold:
\begin{equation}
    \begin{cases}
        \sum_{i'=0}^{m}N_{i'}\theta_{i+i'}=\alpha_i&\forall i=1,\dots, l\\
        \sum_{j'=0}^{m}N_{j'}\phi_{j+j'}=\beta_j&\forall j=1,\dots, m+3-l
    \end{cases}.
\end{equation}
Thus, we can reinterpret (\ref{eqn:EachCase}) in the following manner
\begin{equation}
\label{eqn:muMatrix}
\begin{pmatrix}
\theta_1&\theta_2&\dots&\theta_{m+1}\\
\theta_2&\theta_3&\dots&\theta_{m+2}\\
\vdots&\dots&\ddots&\dots\\
\theta_l&\theta_{l+1}&\dots&\theta_{m+l}\\
\phi_1&\phi_2&\dots&\phi_{m+1}\\
\phi_2&\phi_3&\dots&\phi_{m+2}\\
\vdots&\dots&\ddots&\dots\\
\phi_{m+3-l}&\phi_{m+4-l}&\dots&\phi_{2m+4-l}
\end{pmatrix}\begin{pmatrix}
N_0\\
N_1\\
\vdots\\
N_m\end{pmatrix}=\begin{pmatrix}
\alpha_1\\
\vdots\\
\alpha_l\\
\beta_1\\
\vdots\\
\beta_{m+3-l}
\end{pmatrix}.
\end{equation}
Notice that the system of equations (\ref{eqn:muMatrix}) is a system of equations over $\mathbb{F}_q$ and due to the bound on the denominators of $\theta,\phi,\alpha,\beta$, there are finitely many such systems of equations. This enables us to use a computer program to check if $\mu(\Lambda_{\theta,\phi})\leq q^{-4}$ for rational $\theta,\phi$. In particular, our computations show that in contrast to the real case, in which it is believed that the upper bound of $\mu$ is the unique value to $A\mathbb{Z}^d$, there are many lattices $\Lambda\notin A\mathcal{R}^3$ satisfying $\mu(\Lambda)=q^{-3}$.

Some examples that we found for $q=3$ and $\operatorname{denom}(\theta),\operatorname{denom}(\phi)\leq 10$ include the following:
$$\mu\left(\Lambda_{\frac{1}{x}+\frac{1}{x^2}+\frac{1}{x^4},\sum_{i=1}^{10}\frac{1}{x^i}}\right)=3^{-3},$$
$$\mu\left(\Lambda_{\frac{2}{x^3}+\frac{2}{x^6}+\frac{1}{x^7}+\frac{2}{x^8}+\frac{1}{x^9}+\frac{2}{x^{10}},\frac{1}{x}+\frac{1}{x^2}+\frac{1}{x^3}+\frac{1}{x^5}+\frac{2}{x^6}+\frac{1}{x^9}+\frac{2}{x^{10}}}\right)=3^{-3}.$$
Some examples that we found for $q=2$ and $\operatorname{denom}(\theta),\operatorname{denom}(\phi)\leq 10$ include the following:
$$\mu\left(\Lambda_{\sum_{i=2}^{10}\frac{1}{x^i}-\frac{1}{x^7}}\right)=\mu\left(\Lambda_{\frac{1}{x}+\frac{1}{x^4}+\frac{1}{x^5}}\right)=\mu\left(\Lambda_{\frac{1}{x^9},\frac{1}{x}+\frac{1}{x^3}+\frac{1}{x^7}+\frac{1}{x^8}+\frac{1}{x^10}}\right)=2^{-3}.$$
These computations are evidence that $\mathcal{S}_3$ may be trivial, and we hope to find mathematical tools to prove or disprove Conjecture \ref{minktriv}.
\section{The Mordell Constant}
In this section, we shall prove Theorem \ref{GM}. Due to Theorem \ref{AxcapWR}, Lemma \ref{WRForm}, and $A$ invariance of $\kappa$, for every $\Lambda\in \mathcal{L}_d$, the box $\mathcal{B}=B(0,q^{-1})^d$ is admissible for every $\Lambda$. Thus, $\kappa(\Lambda)\geq q^{-d}$. Hence it suffices to prove that for generic $\Lambda$, $\kappa(\Lambda)\leq q^{-d}$. This follows from Theorem \ref{ClaMCBT}, but we shall provide an alternative proof that $\kappa(\Lambda)\leq q^{-d}$ by using dynamics and improving Dirichlet's theorem. 

We shall use an improved version of Dirichlet's theorem to show that $\kappa(\Lambda)\leq q^{-d}$ for every $\Lambda\in \mathcal{L}_d$. Notice that $A$ invariance and lower semi-continuity of $\kappa$ imply that if $\Lambda$ has a dense $A$ orbit, then, $\kappa(\Lambda)$ must obtain the generic value which is also $\max\mathbf{MG}_d$. 

For $d\geq 2$, denote $n=d-1$. For $\boldsymbol{\theta}=(\theta_1,\dots, \theta_n)\in (x^{-1}\mathcal{O})^n$, let $\Lambda_{\boldsymbol{\theta}}=g_{\boldsymbol{\theta}}\mathcal{R}^{d}$ where 
$$g_{\boldsymbol{\theta}}=\begin{pmatrix}
1&\theta_1&\dots &\theta_n\\
0&1&\dots&0\\
\vdots&\dots&\ddots&\dots\\
0&\dots&\dots&1
\end{pmatrix}.$$
Then, due to a slight variation of Lemma \ref{DenseUni}, for almost every $\boldsymbol{\theta}\in (x^{-1}\mathcal{O})^n$, $\Lambda_{\boldsymbol{\theta}}$ has a dense $A$ orbit, and hence, it suffices to calculate $\kappa(\Lambda_{\boldsymbol{\theta}})$ for such $\boldsymbol{\theta}$.

We shall follow the proof of \cite[Theorem 2.2]{GaGh17} to obtain a better bound on Dirichlet's Theorem. This proof is identical to the proof of \cite[Theorem 2.2]{GaGh17}, and we include it to show how the tighter bound on Dirichlet's Theorem is obtained. This, along with dynamical considerations, will enable us to calculate $\kappa(\Lambda_{\boldsymbol{\theta}})$. 
\begin{theorem}[Improved Dirichlet's Theorem]
\label{improvDirichlet}
Let $\theta_1,\dots, \theta_n\in (x^{-1}\mathcal{O})^n$ and let $t_1,\dots, t_n\geq 0$ be integers. Then, there exist $a,b_1,\dots b_n\in \mathcal{R}$, where $(b_1,\dots b_n)\neq (0,\dots, 0)$, such that
\begin{equation}
    \left|b_1\theta_1+\dots +b_n\theta_n-a\right|\leq q^{-\left(\sum_{i=1}^nt_i+n\right)},
\end{equation}
\begin{equation}
    \vert b_i\vert\leq q^{t_i}.
\end{equation}
\end{theorem}
\begin{proof}
We shall use the observations in \cite{GaGh17}. Let $\alpha=\sum_{i=1}^{\infty}\alpha_ix^{-i}$, $\beta=\sum_{i=1}^{\infty}\beta_ix^{-i}$, and $a=\sum_{i=0}^ka_ix^i,b=\sum_{j=0}^lb_jx^j\in \mathcal{R}$. Then, the coefficient of $x^{-s}$ in $a\alpha+b\beta$ is
$$\sum_{i=0}^ka_i\alpha_{s+i}+\sum_{j=0}^lb_j\beta_{s+j}.$$
Hence, $\vert \langle a\alpha+b\beta\rangle\vert\leq q^{-(m+1)}$ if and only if for each $s=1,\dots, m$, 
$$\sum_{i=0}^ka_i\alpha_{s+i}+\sum_{j=0}^lb_j\beta_{s+j}=0.$$ 
Denote $$\mathbf{a}=\begin{pmatrix}
a_0\\
a_1\\
\vdots\\
a_k
\end{pmatrix}\in \mathcal{R}^{n+1}, \mathbf{b}=\begin{pmatrix}
b_0\\
b_1\\
\vdots\\
b_l
\end{pmatrix}.$$ 
Then, $\vert \langle a\alpha+b\beta\rangle\vert<q^{-m}$ if and only if $\begin{pmatrix}
\mathbf{a}\\
\mathbf{b}
\end{pmatrix}$ is a solution to the linear equation $\begin{pmatrix}
A&B\end{pmatrix}\begin{pmatrix}\mathbf{x}\\
\mathbf{y}
\end{pmatrix}=0$, where
$$A=\begin{pmatrix}
\alpha_1&\dots &\dots &\alpha_{k+1}\\
\alpha_2&\dots &\dots &\alpha_{k+2}\\
\vdots&\dots &\ddots&\dots\\
\alpha_m&\dots &\dots&\alpha_{m+k}
\end{pmatrix}, B=\begin{pmatrix}
\beta_1&\dots &\dots &\beta_{l+1}\\
\beta_2&\ddots &\dots &\beta_{l+2}\\
\vdots&\dots &\ddots&\dots\\
\beta_m&\dots &\dots&\beta_{m+l}
\end{pmatrix}.$$
Returning to the notation of the theorem, take $\vert b_i\vert\leq q^{t_i}$ so that $b_i=\sum_{j=0}^{t_i}b_{i,j}x^j$ corresponds to a vector $\mathbf{b}_i=\begin{pmatrix}
b_{i,0}\\
b_{i,1}\\
\vdots\\
b_{i,t_i}
\end{pmatrix}\in \mathbb{F}_q$. Write $\theta_i=\sum_{j=1}^{\infty}\theta_{i,j}x^{-j}$. Then, the coefficient of $x^{-s}$ in $\sum_{i=1}^nb_i\theta_i$ is equal to
$$\sum_{i=1}^n\sum_{j=0}^{t_i}b_{i,j}\theta_{i,s+j}.$$
Denote $\Theta_{i,s}=(\theta_{i,s},\dots, \theta_{i,s+t_i}).$ Then, $\vert \langle b_1\theta_1+\dots +b_n\theta_n\rangle\vert\leq q^{-(m+1)}$ where $\vert b_i\vert\leq q^{t_i}$ corresponds to the following system of equations
\begin{equation}
    \Theta\mathbf{b}=\begin{pmatrix}
    \Theta_{1,1}&\dots&\dots&\Theta_{n,1}\\
    \Theta_{1,2}&\dots&\dots&\Theta_{n,2}\\
    \vdots&\dots&\ddots&\dots\\
    \Theta_{1,m}&\dots&\dots&\Theta_{n,m}
    \end{pmatrix}\begin{pmatrix}
    \mathbf{b_1}\\
    \vdots\\
    \mathbf{b}_n
    \end{pmatrix}=\boldsymbol{0}.
\end{equation}
Notice that every $\Theta_{i,s}$ has $t_i+1$ columns and hence, $\Theta$ has $\sum_{i=1}^nt_i+n$ columns. Hence, if $m\leq \sum_{i=1}^nt_i+n-1$, then, $\Theta\mathbf{b}=\boldsymbol{0}$ has a nontrivial solution. In particular, there exists a nontrivial solution for $m=\sum_{i=1}^nt_i+n-1$. Hence, there exist $b_i\in \mathcal{R}$ which are not all $0$ such that 
\begin{equation}
    \left|\langle b_1\theta_1+\dots +b_n\theta_n\rangle\right|\leq q^{-\left(\sum_{i=1}^nt_i+n\right)},
\end{equation}
\begin{equation}
    \forall i=1,\dots n, \vert b_i\vert\leq q^{t_i}.
\end{equation}
\end{proof}
We shall now show that for generic $\boldsymbol{\theta}$,  $\kappa(\Lambda_{\boldsymbol{\theta}})=q^{-d}$ by showing that the box
$$\mathcal{B}_0=B(0,q^{-n})\times B(0,1)\times \dots \times B(0,1)$$ 
is not admissible for any unimodular lattice. 
\begin{theorem}
\label{genVal}
For generic $\boldsymbol{\theta}$, $\kappa(\Lambda_{\boldsymbol{\theta}})\leq q^{-d}$.
\end{theorem}
\begin{proof}[Proof of Theorem \ref{GM}]
By Theorem \ref{genVal} and lower semicontinuity of $\kappa$, $\max_{\Lambda\in \mathcal{L}_d}\kappa(\Lambda)\leq q^{-d}$. On the other hand, due to Theorem \ref{AxcapWR} and Lemma \ref{WRForm}, for every $\Lambda\in \mathcal{L}_d$, there exists some $\mathbf{a}\in A$, such that the box $B(0,q^{-1})^d$ is admissible for $\mathbf{a}\Lambda$. Hence, for every $\Lambda\in \mathcal{L}_d$, $\kappa(\Lambda)\geq \operatorname{Vol}(B(0,q^{-1})^d)=q^{-d}$, so that $\kappa(\Lambda_{\boldsymbol{\theta}})=q^{-d}$. Since $\Lambda_{\boldsymbol{\theta}}$ has a dense $A$ orbit, then lower semicontinuity of $\kappa$ implies that $\max\mathbf{MG}_d=\kappa(\Lambda_{\boldsymbol{\theta}})=q^{-d}$. On the other hand, $\min\mathbf{MG}_d\geq q^{-d}$, and thus $\mathbf{MG}_d=\big\{q^{-d}\big\}$. 
\end{proof}
\begin{proof}[Proof of Theorem \ref{genVal}]
We shall first reinterpret Theorem \ref{improvDirichlet} in a different language. Let 
$$A_{-}=\Bigg\{\mathbf{t}=\operatorname{diag}\left(\prod_{i=1}^nt_i^{-1},t_1,\dots, t_n\right):\forall i, \vert t_i\vert<1\Bigg\}<A.$$
Then, Theorem \ref{improvDirichlet} is equivalent to the following: For every $\mathbf{t}\in A_{-}$ and for every $\boldsymbol{\theta}\in (x^{-1}\mathcal{O})^n$, $\mathbf{t}\Lambda_{\boldsymbol{\theta}}\cap \mathcal{B}_0\neq \{0\}$. For $i=1,\dots, n$, let
$$A_{-}^{(i)}=\Bigg\{\mathbf{t}=\operatorname{diag}\left(\prod_{i=1}^nt_i^{-1},t_1,\dots, t_n\right):t_i=1, \forall i\neq j, \vert t_j\vert\leq 1\Bigg\}.$$
Let $\{g_{\mathbf{t}}\}_{\mathbf{t}\in \tilde{\mathcal{K}}}\subseteq A_-\setminus \cup_{i=1}^dA_-^{(i)}$ be an ergodic one parameter flow contained in the interior of the cone $A_-$. Then we can assume that $\boldsymbol{\theta}$ was chosen so the orbit $\{g_{\mathbf{t}}\Lambda_{\boldsymbol{\theta}}\}_{t\in \tilde{\mathcal{K}}}$ is dense. Then, for every $\Lambda\in \mathcal{L}_d$, there exists a sequence $g_{\mathbf{t}_k}\in A_-$ with $d(g_{\mathbf{t}_k},A_-^{(i)})\rightarrow 0$ for every $i=1,\dots, d$, such that $g_{\mathbf{t}_k}\Lambda_{\boldsymbol{\theta}}\rightarrow \Lambda$.
 
For compact $\mathcal{C}\subseteq \tilde{\mathcal{K}}^d$ and $l\geq 1$ define the sets
$$\mathcal{L}_d(\mathcal{C},l)=\{\Lambda\in \mathcal{L}_d:\vert \Lambda\cap\mathcal{C}\vert\geq l\}.$$
Then, by lower semicontinuity, $\mathcal{L}_d(\mathcal{C},l)$ are closed sets. Notice that Theorem \ref{improvDirichlet} implies that for every $\mathbf{t}\in A_{-}$, $\mathbf{t}\Lambda_{\boldsymbol{\theta}}\in \mathcal{L}_d(\mathcal{B}_0,2)$. Hence, if $\Lambda=\lim_{k\rightarrow \infty}g_{\mathbf{t}_k}\Lambda_{\boldsymbol{\theta}}$, then, $\Lambda\in \mathcal{L}_d(\mathcal{B}_0,2)$. Since $\{g_{\mathbf{t}_k}\Lambda_{\boldsymbol{\theta}}\}_{\mathbf{t}_k\in \tilde{\mathcal{K}}}$ is dense, then for every $\Lambda\in \mathcal{L}_d$, $\Lambda\in \mathcal{L}_d(\mathcal{B}_0,2)$. Therefore, the box $\mathcal{B}_0$ is not admissible for any lattice. 
Assume that $\mathcal{B}$ is a box around the origin with $\operatorname{Vol}(\mathcal{B})=q^{-n}$, which is admissible for $\Lambda_{\boldsymbol{\theta}}$. Then, there exists some $\mathbf{a}\in A$ such that $\mathbf{a}\mathcal{B}=\mathcal{B}_0$. Therefore, $\mathcal{B}_0$ is admissible for the lattice $\mathbf{a}\Lambda_{\boldsymbol{\theta}}$, which is a contradiction. Hence, there are no boxes of volume $q^{-n}$ which are admissible for $\Lambda_{\boldsymbol{\theta}}$, and thus, $\kappa(\Lambda_{\boldsymbol{\theta}})\leq q^{-d}$. 
\end{proof}
\bibliography{Ref}
\bibliographystyle{amsalpha}
\end{document}